\documentclass[11pt]{amsart}
\usepackage{amssymb}
\usepackage{graphicx}
\input xy
\xyoption{all}
\theoremstyle{plain}
\newtheorem{theorem}{Theorem}[section]
\newtheorem{corollary}[theorem]{Corollary}
\newtheorem{lemma}[theorem]{Lemma}
\newtheorem{subclaim}[theorem]{Subclaim}
\newtheorem{conjecture}[theorem]{Conjecture}
\newtheorem{proposition}[theorem]{Proposition}
\newtheorem{fact}[theorem]{Fact}
\newtheorem{claim}[theorem]{Claim}
\newtheorem{question}[theorem]{Question}
\theoremstyle{definition}
\newtheorem{definition}[theorem]{Definition}
\newtheorem{remark}[theorem]{Remark}

\newtheorem{assumption}[theorem]{Assumption}
\newtheorem{notation}[theorem]{Notation}
\newtheorem{example}[theorem]{Example}
\newtheorem{examples}[theorem]{Examples}
\theoremstyle{remark}

\newcommand{\id}{\operatorname{id}}
\newcommand{\LS}{\operatorname{LS}}

\newcommand{\cf}{\operatorname{cf}}

\newcommand{\dom}{\operatorname{dom}}
\newcommand{\Aut}{\operatorname{Aut}}

\newcommand{\ftp}{\operatorname{tp}}
\newcommand{\gatp}{\operatorname{ga-tp}}
\newcommand{\tp}{\operatorname{ga-tp}}
\newcommand{\gaS}{\operatorname{ga-S}}
\newcommand{\Hanf}{\operatorname{Hanf}}

\newcommand{\conc}{\Hat{\ }}

\newcommand{\sq}[2]{\sideset{^{#1}}{}{\operatorname{#2}}}

\newcommand{\Mod}{\operatorname{Mod}}

\newcommand{\Pm}{\mathcal{P}^{-}}

\renewcommand{\phi}{\varphi}

\newcommand{\Union}{\bigcup}
\newcommand{\initial}\lessdot
 \newcommand{\infinity}{\infty}
\newcommand{\K}{\operatorname{\mathcal{K}}}
\def\l{\langle}
\def\r{\rangle}

\newcommand{\C}{\mathfrak C}

\def\?{?\vadjust
{\vbox to 0pt{\vskip-7pt\hbox to 1.1\hsize{\hfill\huge ?!}}}}

\begin{document}

\title{Galois-stability for Tame Abstract Elementary Classes}

\author{Rami Grossberg}
\email[Rami Grossberg]{rami@andrew.cmu.edu}
\address{Department of Mathematics\\
Carnegie Mellon University\\
Pittsburgh PA 15213}

\author{Monica VanDieren}
\email[Monica VanDieren]{mvd@umich.edu}
\address{Department of Mathematics\\
University of Michigan\\
Ann Arbor MI 48109-1109}


\date{September 22, 2005}
\thanks{This paper is part of the second author's Ph.D. thesis written under
the guidance of Rami Grossberg}

\maketitle

\begin{abstract}

We introduce \emph{tame abstract elementary classes} as a generalization
of all cases of abstract elementary classes that are known to permit
development of stability-like theory.  In this paper we  explore stability
results in this new context.  We assume that
$\K$ is a tame abstract elementary class satisfying the amalgamation
property with no maximal model. The main results include:

\begin{theorem}\label{bounded non-split}
Suppose that $\K$ is not only tame, but $\LS(\K)$-tame.
If
$\mu\geq \Hanf(\K)$ and
$\K$ is Galois
stable in $\mu$, then $\kappa_\mu(\K)<\beth_{(2^{\Hanf(\K)})^+}$, where
$\kappa_{\mu}(\K)$ is a  relative of $\kappa(T)$ from first order logic.
\index{$\kappa^\beta_\mu(\K)$}

$\Hanf(\K)$ is the Hanf number of the class $\K$. It is known that
$\Hanf(\K)\leq\beth_{({2^{\LS(\K)}})^+}$
\end{theorem}

The theorem generalizes a result from \cite{Sh3}.  It is used to prove
both the existence of Morley sequences for non-splitting (improving Claim
4.15 of \cite{Sh 394} and a result from \cite{GrLe1}) and
%
%\begin{theorem}\label{mseq thm}
%There exists a cardinal $\mu_0(\K)$ such that for every
% $\mu\geq \mu_0(\K)$ and every $M\in\K_{>\mu}$, $A,I\subset M$
%such that
%$|I|\geq\mu^+>|A|$, if $\K$ is Galois-stable in $\mu$,
%then there exists $J\subset I$ of cardinality $\mu^+$,
%Galois-indiscernible sequence over $A$.  Moreover $J$ can be chosen
%to be a Morley sequence  over $A$.\index{Galois-indiscernible
%sequence!existence}\index{Morley sequence!existence}
%\end{theorem}
%
%This result strengthens
%Claim 4.16 of \cite{Sh 394} as we do not assume categoricity.  This is
%also an improvement of  a result from
%\cite{GrLe1} concerning the existence of indiscernible sequences.
the following initial step towards a stability spectrum theorem for
tame classes:

\begin{theorem}
If $\K$ is Galois-stable in some $\mu>\beth_{(2^{\Hanf(\K)})^+}$, then
$\K$ is stable in every $\kappa$ with $\kappa^\mu=\kappa$.
E.g. under \emph{GCH} we have that $\K$ Galois-stable in $\mu$ implies that
$\K$ is Galois-stable in $\mu^{+n}$ for all $n<\omega$.
\end{theorem}

%We also improve a result from \cite{ShVi 635} by removing the assumption
%of categoricity and GCH:
%
%\begin{theorem}
%Let $\K$ be an abstract elementary class without maximal models and
%suppose it is
%Galois-stable in $\mu$.  If $\K$ has the disjoint amalgamation property
%then for every
%$M\in\K_\mu$ there exists
%$M^*\succneqq_{\K}M$, universal over $M$ of cardinality $\mu$.
%\end{theorem}

\end{abstract}

%%%%%%%%%%%%%%%%%%%%%%%%%%%%%%%
\bigskip
\section*{Introduction} \label{s:introductionm}

In the last twenty years most of Shelah's effort in model theory was
in developing classification theory for abstract elementary classes (Some important  mile
stones are:
\cite{Sh 87a},
\cite{Sh 87b},\cite{Sh 88},
\cite{Sh 300},
\cite{Sh 394},
\cite{Sh 576},
\cite{Sh 600}
and
\cite{Sh705}).  
An abstract elementary class (AEC) is a
class of structures of the
same similarity type endowed with a morphism satisfying natural properties
such as closure under directed limits (see \cite{Gr1} for an
introduction).  The notion of an AEC (introduced by Shelah in \cite{Sh 88}) is broad enough to
capture classes axiomatized by 
non-first-order  logics including $L_{\omega_1,\omega}(\mathbf{Q})$.
  To date, there is no known stability theory or even a 
categoricity theorem for abstract elementary classes without some additional strong
assumptions.

The most general of
Shelah's attempts to develop a stability theory is
\cite{Sh 600} where Shelah works in the context of
AECs with the amalgamation and joint embedding properties and a very nicely 
behaved forking-like relation.  He
introduces the notion of good-frame which is a generalization of
a first-order superstable theory  (the definition alone needs a couple of
pages).  Most of his 100-page paper is dedicated to the derivation of
structural  results  under the assumption that the AEC has a forking-like
relation.

In our paper we deal with a much wider class of AECs than Shelah in
\cite{Sh 600} as we do not assume as part of the setting built-in
forking and regular types.  We deal with stable but tame
classes (see Definitions \ref{stable defn} and \ref{tame defn}), 
in other words, stable classes where inequality of types behaves
locally.  
Tameness is a general property capturing finite 
diagrams (homogeneous model theory).
While it is not difficult to see that excellent classes of
atomic models of a first-order theory (from \cite{Sh 87a} and
\cite{Sh 87b}) are tame, Grossberg and Alexei Kolesnikov show
that
  every abstract excellent class is also tame (see
\cite{GrKo}).
%Thus our context is a strict generalization of
%\cite{Sh705} where the work is done on AECs
%which are also good frames and
%excellent.
%%
%
Very recently using ideas from  \cite{GrKo},  Andr\'{e}s 
Villaveces and Pedro 
Zambrano  discovered that  the class of structures obtained from 
Hrushovski's fusion of strongly minimal first-order theories is another 
source of examples of tame AECs.  See \cite{ViZa} and \cite{Ba}.  

Despite the generality of our context, recent progress on a Morley's 
Theorem for tame abstract elementary classes provides evidence
that
a stability-theory can be developed.
A short time after submitting this paper for publication we
realized that 
using the splitting machinery developed here,
we can prove in ZFC an instance of Shelah's categoricity
conjecture for tame classes (see \cite{GrVa1} and \cite{GrVa2}).

%For many years the best approximation to Shelah's categoricity conjecture 
%was the  following beautiful theorem from \cite{Sh285}:
%%
%
%\begin{theorem}[Makkai and Shelah 1990]
%      Let $\K$ be an AEC,  \emph{  $\kappa$ a strongly compact cardinal}
%such that $\LS(\mathcal{K})<\kappa$.  Let $\mu_0:=\beth_{(2^\kappa)^+}$.
%If
%$\K$ is categorical in some
%      $\lambda^{+}>\mu_0$ then $\K$ is categorical in every $\mu\geq\mu_0$.
%\end{theorem}
%
%  Proposition 1.13 of \cite{Sh285} asserts  (using the assumption that 
%$\kappa$ is strongly compact)
%that any AEC $\K$ as above has the amalgamation property (for models of 
%cardinality
%$\ge\kappa$). Since Galois types in this context are sets of 
%$L_{\kappa,\kappa}$ formulas the class is trivially $\kappa$-tame.

%      In
%\cite{GrVa1} Grossberg and VanDieren proved (in ZFC)  a case of  Shelah's
%categoricity conjecture for  tame AECs with  the amalgamation property
%which implies the conclusion of the Makkai and Shelah theorem. Thus the
%tameness assumption enables upward categoricity argument  in ZFC (instead 
%of the
%large cardinal assumption).  This is also an (upward)  extension of Shelah's
%main theorem from \cite{Sh 394}.
%%
%
\begin{theorem}[Grossberg and VanDieren 2003]
     Let $\K$ be an AEC, $\kappa:=\beth_{(2^{\LS(\K)})^+}$.  Denote by
$\mu_0:=\beth_{(2^\kappa)^+}$.  Suppose that $\K_{>\kappa}$ has the
amalgamation property and is $\kappa$-tame.
If
$\K$ is categorical in some
      $\lambda^{+}>\mu_0$ then $\K$ is categorical in every $\mu\geq\mu_0$.
\end{theorem}

This work removes the set-theoretic assumption of a strongly compact 
cardinal from the work of Makkai and Shelah \cite{Sh285}.

In addition to examining splitting in this paper, 
we prove the existence of Morley sequences from the 
assumption of stability in tame AECs.
 We plan to use
Morley sequences for non-splitting to define a
dividing-like concept for these classes and to prove a stability
spectrum theorem for AECs.  More work on the stability spectrum theorem 
(extension of our results under the additional assumption of 
$\LS(\K)=\aleph_0$)
appears in \cite{BaKuVa}.

We are grateful to John Baldwin for reading several  preliminary
versions and suggesting us to consider the stability spectrum problem for
AECs.  The referee's comments and questions led to significant improvement
of the presentation.

%%%%%%%%%%%%%%%%%%%%%%%%%%%%%%%
\bigskip
\section{Background} \label{s:backgroundm}

Much of the necessary background for this paper can be found in
\cite{Gr1}.  We will use $\alpha,\beta,\gamma,i,j$ to denote ordinals and
$\kappa,\lambda,\mu,\chi$ will be used for cardinals.
 We will use $(\K,\prec_{\K})$ to denote an abstract
elementary class and $\K_\mu$ is the subclass of models in $\K$ of 
cardinality
$\mu$.  Similarly we define $\K_{<\mu}$ and $\K_{\geq\mu}$.  Models are
denoted by
$M,N$ and may be decorated with superscripts and subscripts.  Sequences of
elements from
$M$ are written as
$\bar a,\bar b,\bar c,\bar d$ and we write $\sq\beta M$ to denote the
collection of all sequences of length $\beta$ made up of elements from
$M$. The letters $e,f,g,h$ are reserved for $\K$-mappings and $\id$ is the
identity mapping.

For the remainder of this paper we will fix $(\K,\prec_{\K})$ to be an
abstract elementary class.

We begin by reviewing the definition of
Galois-type, since we will be considering variations of the underlying
equivalence relation
$E$ later in this paper.

\begin {definition}\label{E}\index{$E$, binary relation}
Let $\beta>0$ be an ordinal.
For triples $(\bar a_\ell, M, N_l)$ where $\bar a_\ell\in \sq{\beta}N_\ell$ 
and
$M_\ell\prec_{\K}N_\ell\in\K$ for
$\ell=0,1$,
we define a binary relation $E$ as follows:
$(\bar a_0, M, N_0)E(\bar a_1, M, N_1)$ iff
 and there exists $N\in\K$ and $\K$-mappings
$f_0,f_1$ such that
$f_\ell:N_\ell\rightarrow N,$    $f_\ell\restriction M=\id_M$ for 
$\ell=0,1$ and
$f_0(\bar a_0)=f_1(\bar a_1)$:

\[
\xymatrix{\ar @{} [dr] N_1
\ar[r]_{f_1}  &N \\
M \ar[u]^{\id} \ar[r]_{\id}
& N_0 \ar[u]_{f_0}
}
\]

\end{definition}
 The relation $E$ is used to define Galois-types in classes that satisfy
the amalgamation property.
\begin{definition}
\begin{enumerate}
\item Let $\mu\geq \LS(\K)$.  We say that $\K$ has the\\
 \emph{
$\mu$-amalgamation property} iff for any $M_\ell\in \K_\mu$ (for $\ell\in
\{0,1,2\}$) such that $M_0\prec_{\K} M_1$ and $M_0\prec_{\K}M_2$ there are
$N\in\K_\mu$ and $\K$-embeddings $f_\ell:M_\ell
\rightarrow N$ such that $f_\ell\restriction M_0=\id_{M_0}$
for $\ell =1,2$.

\item A model $M_0\in \K_\mu$ satisfying the above requirement is called
an
\emph{amalgamation base}.

\item $\K$ has the \emph{amalgamation property} iff $\K$ has the
$\mu$-amalgamation property for all $\mu\geq \LS(\K)$.
\end{enumerate}
\end{definition}

\begin{remark}
$E$ is an equivalence relation on the class of triples of the form
$(\bar a, M, N)$ where $M\prec_{\K}N$, $\bar a\in N$ and
both $M$ and $N$ are
amalgamation bases.  When $N$ is not an amalgamation base,
$E$ may fail to be transitive, but the transitive closure of $E$
could be used instead.

\end{remark}

\begin{remark}
While, the focus of this paper is on classes with the amalgamation
property, several of the proofs in this paper can be adjusted to the
context of abstract elementary classes with density of amalgamation bases
as in
\cite{ShVi 635} and \cite{Va}.
\end{remark}

We will make the following assumption for the remainder of the paper:

\begin{assumption}

$\K$ satisfies the amalgamation property.
\end{assumption}
Under the assumption of the amalgamation property,
there exists a large model-homogeneous model, which we will denote by $\C$
and call the monster model.
A \emph{$\mu$-model homogeneous model} is a model $M$ in which for every
$N\prec_{\K}M$ of cardinality $\mu$ and every extension $N_1$ of
cardinality $\mu$, there exists a $\K$-embedding $f:N_1\rightarrow M$
with $f\restriction N=\id_N$.
$M$ is model homogeneous when it is $\mu$-model homogeneous for every
$\mu<\|\C\|$.  All models and sequences of elements will be
assumed to come from
$\C$. Thus,
$E$ is an equivalence relation and we can now define types as equivalence
classes of $E$.

In order to avoid confusing this non-elementary notion of ``type'' with
the conventional first-order
one (i.e. set of formulas) we will follow
\cite{Gr1} and \cite{Gr2} and
 introduce it below under
the name of Galois-type.

\begin{definition}\label{defn of types}\index{Galois-type}\index{type,
Galois} Let $\beta$ be a positive ordinal (can be one).
\begin{enumerate}
\item For $M,N\in\K$ and $\bar a\in \sq{\beta}N.$
The \emph{Galois-type of $\bar a$ in $N$ over $M$}, written
$\tp(\bar a/M,N)$, is defined to be $(\bar a,M,N)/E$.
\item We abbreviate $\tp(\bar a/M,N)$ by $\tp(\bar a/M)$.\index{$\tp(\bar
a/M)$}
\item For $M\in\K$,
$$\gaS^\beta(M):=\{\tp(\bar a/M,N)\mid M\prec
N\in\K_{\|M\|},
\bar a\in \sq{\beta}N\}.$$\index{$\gaS^\beta(M)$}
  We write $\gaS(M)$ for
$\gaS^1(M)$.\index{$\gaS^1(M)$}\index{$\gaS^\beta(M)$}
\item Let $p:=\tp(\bar a/M',N)$ for $M\prec_{\K}M'$ we denote by
$p\restriction M$ the type
$\tp(\bar a/M,N)$.  The \emph{domain of $p$} is denoted by $\dom (p)$ and
it is by definition $M'$.
\index{Galois-type!restriction}\index{Galois-type!domain}
\item  Let $p=\tp(\bar a/M,N)$, suppose that $M\prec_{\K}N'\prec_{\K}N$
and let $\bar b\in \sq{\beta}N'$ we say that \emph{$\bar b$ realizes $p$}
iff
$\tp(\bar b/M,N')=p\restriction M$.\index{Galois-type!realized}
\item For types $p$ and $q$, we write $p\leq q$ if $\dom(p)\subseteq
\dom(q)$ and there exists $\bar a$ realizing $p$ in some $N$ extending
$\dom(p)$ such that $(\bar a,\dom(p), N) = q\restriction
\dom(p)$.\index{Galois-type!extension}
\end{enumerate}
\end{definition}

\begin{definition}\label{stable defn}
We say that $\K$ is
\emph{$\beta$-Galois-stable in $\mu$} \index{$\beta$-Galois-stable in
$\mu$}\index{Galois stable!$\beta$-stable in $\mu$}if for  every
$M\in\K^{}_\mu$,
$|\gaS^\beta(M)|=\mu$.  The class $\K$ is \emph{Galois stable in $\mu$}
iff
$\K$ is $1$-stable in $\mu$.
\end{definition}

\begin{remark}\label{otop remark}
Let $T$ be a stable, countable, first-order, complete theory.  Set
$\K:=\Mod(T)$ and $\prec_{\K}$ the usual elementary submodel relation.
Take $\mu=2^{\aleph_0}$.  While $\K$ is $1$-Galois-stable in $\mu$ (in fact 
it
is $n$-Galois-stable in $\mu$ for all $n<\omega$), it can be shown that
$\K$ is $\omega$-Galois-stable in $\mu$ iff $\K$ has the ndop and notop.

\end{remark}

While there is a nice relationship between
$\beta$-stability and $\beta'$-stability for first order theories, the
corresponding
relationship between
$\beta$-stability and $\beta'$-stability for
AECs in general is unknown and probably fails.

\begin{notation}
Although not standard in AECs with the amalgamation property, we will use
the notation $\gatp(\bar a/\emptyset)$ to denote the orbit of $\bar a$
under $\Aut(\C)$.  This notation is only used as a device in the proof of
Theorem \ref{tame and stable imply kappa exists}.
\end{notation}
%%%%%%%%%%%%%%%%%%%%%%%%%%%%%%%
\bigskip
\section{Saturated and Limit Models} \label{s:sat}

In this section we recall that limit models (see Definition \ref{defn
limit}) can be used as a substitute for saturated models.  Later, we will
see that limit models will be more convenient in situations
where the first order case would call for saturated models.

\begin{definition}\label{defsat}
We say that $M\in\K_{>\LS(\K)}$ is \emph{Galois-saturated}\index{Galois
saturated} if for every
$N\prec_{\K}M$ of cardinality $<\|M\|$, and every
$p\in \gaS(N)$, we have that $M$ realizes $p$.
\end{definition}

\begin{remark}
When $\K=\Mod(T)$ for a first-order $T$, using the compactness theorem
one can show (Theorem 2.2.3 of \cite{Gr1})
that for  $M\in \K$, the model $M$ is Galois-saturated
 iff  $M$ is saturated in the
first-order sense.
\end{remark}

%A central property of AECs with the amalgamation property is:
%
%\begin{fact}[Shelah \cite{Sh 300}]\label{sat iff mh}  Let
%$\lambda>\LS(\K)$. Suppose that
%$\K$ has the amalgamation property and $N\in \K_\lambda$.
%The following are equivalent
%\begin{enumerate}
%\item $N$ is Galois-saturated.
%\item  $N$ is model-homogenous. \index{model-homogenous}I.e. if
%$M\prec_{\K}N$ and
%$M'\succ M$ of cardinality less than $\lambda$ then
%there exists a $\K$-embedding over $M$ from $M'$ into $N$.
%\end{enumerate}
%\end{fact}
%The originally published argument of \cite{Sh 300} is flawed. A
%complete
%and correct proof of Fact \ref{sat iff mh} can be found in \cite{Sh 576}.
%See also
%\cite{Gr1}.
%
%%\begin{remark}
%
%%\end{remark}

We begin with universal extensions which are used to build limit models.
A universal extension captures some properties of saturated models without
referring explicitly to types. The notion of universality over countable
models was first analyzed by Shelah in Theorem 1.4(3) of \cite{Sh 87a}.

\begin{definition}
\begin{enumerate}
\item \index{universal over!$\kappa$-universal over}
Let $\kappa$ be a cardinal $\geq \LS(\K)$.
We say $N$ is \emph{$\kappa$-universal over $M$} iff
for every $M'\in\K_{\kappa}$ with $M\prec_{\K}M'$ there exists
a $\K$-embedding $g:M'\rightarrow N$ such that
$g\restriction M=\id_{M}$:

\[
\xymatrix{\ar @{} [dr] M'
\ar[dr]^{g}  & \\
M \ar[u]^{\id} \ar[r]_{\id}
& N
}
\]

\item \index{universal over}

We say $N$ is \emph{universal over $M$} or $N$ is \emph{a universal
extension of $M$} whenever
$N$ is $\|M\|$-universal over $M$.

\end{enumerate}
\end{definition}

\begin{remark}
Notice that
the definition of \emph{$N$ universal over $M$} requires all extensions
of $M$ of cardinality $\|M\|$ to be embeddable into $N$.  Variants of this
definition in the literature often involve
$\|M\|<\|N\|$.  We would like to emphasize that in this paper we are
concentrating on the case when
$\|M\|=\|N\|$.
\end{remark}

\begin{examples}\begin{enumerate}
\item

Suppose that $T$ is a first-order complete theory that is stable in some
regular $\mu$. Then every model $M$ of $T$ of cardinality $\mu$ has
an elementary extension $N$ of cardinality $\mu$ which is universal over
$M$. To see this define an elementary-increasing and continuous chain of
models of $T$ of cardinality $\mu$,
$\langle N_i\mid i<\mu\rangle$ such that $N_{i+1}$ realizes all types
over $N_i$.  Let $N=\Union_{i<\mu}N_i$.  By a back-and-forth
construction, one can show that $N$ is universal over $M$.

\item   Let $\K$ be the class of
atomic models of a complete first-order theory in a countable language and
$\prec_{\K}$ the usual notion of first-order elementary submodel.
If $1\leq \mbox{I}(\aleph_1,\K)<2^{\aleph_1}$ and
$2^{\aleph_0}<2^{\aleph_1}$ then for every  $M\in \K_{\aleph_0}$ there
exists a proper extension $N\in \K_{\aleph_0}$  universal over $M$, by
Theorem 1.4(3) of \cite{Sh 87a}.
\end{enumerate}
\end{examples}

\begin{definition}\label{defn limit}
For $M',M\in\K_\mu$ and $\sigma$ a limit ordinal with
$\sigma<\mu^+$,
we say that $M'$ is a \emph{$(\mu,\sigma)$-limit over $M$}\index{limit
model!$(\mu,\sigma)$-limit model over $M$} iff there exists a
$\prec_{\K}$-increasing and continuous sequence of models $\langle M_i\in
\K_{\mu}\mid i<\sigma
\rangle$
such that
\begin{enumerate}
\item $M= M_0$,
\item $M'=\Union_{i<\sigma}M_i$
and
\item\label{univ cond in defn} $M_{i+1}$ is universal over $M_i$.
\end{enumerate}
\end{definition}

\begin{remark}
Notice that in Definition \ref{defn limit}, for $i<\sigma$ and
$i$ a limit ordinal, $M_i$ is a $(\mu,i)$-limit model.

\end{remark}

\begin{definition}\label{limit defn}
We say that $M'$ is a \emph{$(\mu,\sigma)$-limit}\index{limit
model!$(\mu,\sigma)$-limit} iff there is some
$M\in\K$ such that $M'$ is a $(\mu,\sigma)$-limit over $M$.
\end{definition}

When $\K=\Mod(T)$ for a first-order and stable $T$ then automatically (by
Theorem III.3.12 of \cite{Shc}):
\[
M\in \K_\mu\text{ is saturated }\implies
\begin{array}[t]{l}
M \text{ is }(\mu,\sigma)
\text{-limit for all }\sigma<\mu^+\\
\text{of cofinality
}\geq\kappa(T).\end{array}
\]

The uniqueness  of $(\mu,\sigma)$-limit models is an interesting problem:
The statement is
\[
(*)\quad \text{If }N_\ell \text{ are } (\mu,\sigma_\ell)
\text{-limit
models over } M
\text{ then }N_1\cong_M N_2.
\]
A standard back and forth construction gives a uniqueness proof when
$\cf(\sigma_1)=\cf(\sigma_2)$.

Problem $(*)$ for $\sigma_1\neq \sigma_2$ is important.  Let's make
some basic observations:

1.  Suppose $\K=\Mod(T)$ and $T$ is a
complete first-order theory.  If $\mu$ is regular and $\sigma=\mu$ then the
$(\mu,\sigma)$-limit model (if it exists) 
is the saturated model of $T$ of cardinality
$\mu$.

%2.  Let $T:=T_{ord}$. If $N$ is $(\mu,\omega)$-limit over $M$ then $N$ is
%not saturated.  Suppose
%$\{M_n\mid n<\omega\}$ is a chain of universal models such that
%$N=\bigcup_{n<\omega}M_n$ and
%$M_0=M$.   Let $a_n\in M_{n+1}$ realizing the type
%\[
%\{x>a\mid a\in |M_n|\}.
%\]
%Now since the type $\{x>a_n\mid  n<\omega\}$ is omitted in $N$ the model
%is not even
%$\aleph_1$-saturated.

%Essentially the same example shows that for any unstable theory.
%If $\sigma_1<\sigma_2$ and $\sigma_2$ is regular then
%$N_1\not\cong _M N_2$ when $N_\ell$ is $(\mu,\sigma_\ell)$-limit  over
%$M$.

%3.  A similar argument shows that $(\mu,\omega)$-limit models are not
%isomorphic
%$(\mu,\mu)$-limit models for $T$ unsuperstable theories (at least when
%$\mu $ is regular).

%4.  When $T$ is $\aleph_0$-stable then ...  ?? Uniq implies categoricity? 
%probably not.
%
%5.  For f.o. is $T \text{ super stable }\iff (*)$?

2.  Showing that two $(\mu,\sigma)$-limit models are isomorphic for
different
$\sigma$s is a central (dichotomy) property even for first-order
theories: E.g. Let $T$ be a countable (for simplicity), stable,
non-superstable first-order theory. Suppose $T$ is stable in $\mu$.  By
Harnik's theorem \cite{Ha} $T$ has a saturated model of cardinality
$\mu$, in fact for $\sigma<\mu^+$ the theory has $(\mu,\sigma)$-limit
models. Since $T$ is not superstable
one can modify the construction in Thm 13.1 of \cite{AlGr} to construct a $(\mu,\omega)$-limit model which is not saturated.
 there exists a union of an
$\omega$-sequence of saturated models which is not saturated (see
\cite{AlGr}).  The union of such as sequence 
is not isomorphic to a $(\mu,\omega_1)$-limit model
(that was obtained by taking a union of saturated models, see Theorem
III.3.12 of \cite{Shc}).  Thus,
when $T$ is countable and stable but not superstable,
the
saturated model of cardinality $\mu$ is a
$(\mu,\omega_1)$-limit but not a
$(\mu,\omega)$-limit model.
%\?Rami, is the example from \cite{AlGr} a $(\mu,\omega)$-limit?

We will see that under the amalgamation property, the existence of limit
models follows from stability.  Limit models have been applied in several
contexts including \cite{Sh 394},\cite{Sh 576} and \cite{Sh 600}.
Even without the amalgamation property,
the
existence
\cite{ShVi 635} and uniqueness \cite{Va} have been shown for
categorical AECs under some mild set-theoretic assumptions.

Lately Grossberg, VanDieren and Villaveces proved uniqueness of limit
models under mild model-theoretic assumptions on the AEC without assuming
categoricity or working outside  of ZFC (see \cite{GVV}).

We now describe how one proves the existence of universal extensions 
and limit models in stable AECs.
It turns out that for first-order  theories stable in $\mu$ any model of 
cardinality $\mu$ has a
universal  extension of cardinality $\mu$.
The result for AECs is claimed
 in [Sh600] without a proof.

\begin{claim}[Claim 1.15.1 from \cite{Sh 600}]\label{existence of univ}
Suppose $\K$ is an abstract elementary class with the amalgamation property.
If $\K$ is Galois-stable in $\mu$, then for every $M\in\K_\mu$, there
exists
$M'\in\K_\mu$ such that $M'$ is universal over $M$.
Moreover $M'$ can be chosen to be a $(\mu,\sigma)$-limit over $M$ for
any $\sigma<\mu^+$.
\end{claim}

We provide a proof of this claim here for completeness.

As it is easy to prove:

\begin{proposition}
Suppose $\kappa\geq\LS(\K)$, $\K$ has the $\kappa$-amalgamation property
and
$M,N\in\K_\kappa$.
If $M\precneqq_{\K}N$ and $N$ is maximal (i.e. there is no
$N'\succneqq_{\K}N$) then $N$ is universal over $M$.
\end{proposition}

%In this setting it is natural to assume that there are no maximal models.
% under slightly
%stronger hypothesis.

%We provide a proof of the existence of universal extensions in stable AECs
%that satisfy the disjoint amalgamation property.
%
%\begin{definition}
%$\K_\mu$ has the \emph{disjoint amalgamation property}
%\index{disjoint amalgamation property} iff for any $N_\ell\succneqq_{\K}M$ 
%(all of
%cardinality $\mu$) there are $N^*\in \K_\mu$, $N^*\succ_{\K}N_2$ and a 
%$\K$-embedding
%$g:N_1\rightarrow N^*$ such that
%\[
%\xymatrix{N_1 \ar[r]^{g} & N^{*} \\
%M \ar[u]^{\id} \ar[r]_{\id} & N_2 \ar[u]_{\id}
%}
%\]
%commutes and $g[N_1]\cap N_2=M$.
%
%  We denote this property by \emph{$\mu$-DAP}.  \index{$\mu$-DAP}
%\end{definition}
%
%
%
%\begin{remark}
%For first-order complete $T$ using Robinson's consistency lemma $\Mod(T)$ 
%has the amalgamation
%property
%for all $\mu\geq |T|$.  By a theorem  of Lascar and Poizat from
%\cite{LaPo} if $T$ is
%a complete first-order theory then $\Mod(T)$ has the $\mu$-disjoint
%amalgamation property for all
%$\mu\geq |T|$.  See also Theorem 6.4.3 in Hodges \cite{Ho}.
%\end{remark}
%

%Since always we assume that $\K$ has the amalgamation property and the 
%classes we
%consider are always categorical in $\LS(\K)$ (notice that taking the Scott 
%sentence
%of an uncountable large structure does not change the spectrum
%function for uncountable cardinals).  It follows that automatically 
%(mainly because
%of the AP) we may assume:

For the rest of this paper we assume:

\begin{assumption}
$\K$ does not have maximal models.
\end{assumption}

This is a reasonable assumption, since (in the presence of the
disjoint amalgamation property) it holds automatically for categorical
classes or classes with the joint embedding property.
%In
%the special case when
%$\K=\Mod(T)$ for an
%$L_{\omega_1,\omega}$-theory and
%$\prec_{\K}$  is given by a countable fragment then we may replace $\K$ 
%with some
%$\K'$ which is $\Mod(\psi)$ for some Scott sentence (i.e. $\K'$ is
%$\aleph_0$-categorical) of an uncountable element of
%$\K$ and both $\K$ and its subclass $\K'$ have the same spectrum function 
%on
%uncountable cardinals.

\begin{theorem}\label{exist univ}
Let $\K$ be an AEC which is 
$1$-Galois-stable in
$\mu$.  If
$\K$ has the amalgamation property  then for every
$M\in\K_\mu$ there exists
$M^*\succneqq_{\K}M$, universal over $M$ of cardinality $\mu$.
\end{theorem}

\begin{proof}

Define an increasing continuous chain $\langle M_i\mid
i<\mu\rangle\subseteq
\K_\mu$ such that
$M_0:=M$ and every $p\in \gaS(M_i)$ is realized in $M_{i+1}$.  Let 
$M^*:=\bigcup_{i<\mu}M_i$.

We now show that $M^*$ is universal over $M$.
Let $N\succ_{\K}M$ be a given model of cardinality $\mu$, we'll construct an
embedding $f:N\rightarrow M^*$ s.t. $f\restriction M=\id _M$.

Fix $\{a_i\mid i<\mu\}$ an enumeration of $|N|-|M|$.

By induction on $i<\mu$ define $\prec_{\K}$-increasing and continuous
chains\\
$\{N_\ell^i\mid i<\mu,\ell=0,1\}\subseteq \K_\mu$ and $\{f_i\mid i<\mu\}$ 
such that
\begin{enumerate}

\item $N_0^i\prec_{\K}N_1^i$,

\item  $N_0^0=M$, $N_1^0=N$

\item  $a_i\in N_0^{i+1}$

\item  $f_i: N_0^i\rightarrow M_i$.

\end{enumerate}

Clearly this is sufficient since $(\bigcup_{i<\mu}f_i)\restriction N$ is as 
required.
By continuity it is enough
to describe the construction
at successor stages.

If $a_j\in N_0^j$ let $N_\ell^{j+1}:=N_\ell^j$.

Otherwise, denote by $M^j_0\prec_{\K}M_j$ the image of $N_0^j$
under
$f_j$, let
$g\supseteq f_j$ and $\overline{M}_1^j\succ_{\K} M_0^j$ be such that 
$g:N_1^j\cong
\overline{M}_1^j$.

By the amalgamation property there are $M^j_1$ and
$\overline{g}:\overline{M}_1^j\rightarrow M_1^j$ whose restriction is the 
identity over
$M_1^j$ such that
 the diagram
\[
 \xymatrix{\ar @{} [dr] N_1^j
 \ar[r]^{g}  & \overline{M}_1^j \ar[r]^{\overline{g}} & {M}_1^j\\
N_0^j \ar[u]^{\id} \ar[r]_{f_j}
& M_0^j \ar[u]_{\id} \ar[r]_{\id} & M_j\ar[r]_{\id}\ar[u]_{\id} & M
 }
\]
commutes.

Since $a_j\in N_1^j- N_0^j$, $g(a_j)\in
\overline{M}_1^j-M_0^j$.  We can choose
$\bar g$ such that
$\overline{g}(g(a_j))\not\in M_j$.  Let
$p:=\gatp(\overline{g}(g(a_j))/M_j,M_1^j)$. By the construction of
$\langle M_i\rangle_{i<\mu}$ there exists $b\in M_{j+1}$ realizing the
type
$p$.  By definition of Galois-types we have
$(M_j,M_1^j,\overline{g}(g(a_j)))E(M_j,M_{j+1},b)$. Using the
definition of
$E$ (Definition \ref{E}), there are
$N^{**}\in \K_{\mu}$ and mappings $h_1,h_2$ such that
the diagram
\[
 \xymatrix{\ar @{} [dr] M_1^j
 \ar[r]^{h_2}  &N^{**} \\
M_j \ar[u]^{\id} \ar[r]_{\id}
& M_{j+1} \ar[u]_{h_1}
 }
\]
commutes and
\[
(*)\quad h_2(\bar g(g_j(a_j)))=h_1(b).
\]
We may assume that $h_2=\id_{M_1^j}$.
Thus by gluing  the last two diagrams together we get that
the diagram
\[
 \xymatrix{\ar @{} [dr] N_1^j
 \ar[r]^{g\circ \overline{g}}  & M_1^j \ar[r]^{\id} & N^{**}\\
N_0^j \ar[u]^{\id} \ar[r]_{f_j}
& M_j \ar[u]_{\id} \ar[r]_{\id} & M_{j+1}\ar[u]_{h_1}\ar[r]_{\id} & M
 }
\]
commutes.  Now pick $N_1^{j+1}\succ_{\K} N^j_1$ and
$h\supseteq g\circ \overline{g}$ such that
$h:N_1^i\cong N^{**}$.

So we have that
\[
 \xymatrix{\ar @{} [dr] N_1^{i+1}\ar[drr]^{h}\\
N_1^j\ar[u]^{\id}
 \ar[rr]_{g\circ \overline{g}}  &  & N^{**}\\
N_0^j \ar[u]^{\id} \ar[rr]_{f_j}
&  & M_{j+1}\ar[u]_{h_1}
 }
\]
commutes.  Let $N_0^{j+1}:=h^{-1}[h_1[M_{j+1}]]$.\\
Notice that by the above diagram, we get that
 $N_0^{j+1}\supseteq N_0^j$.
Now we are ready to define the mapping $f_i$, let
$f_{j+1}:=h_1^{-1}\circ (h\restriction N_0^j)$. It is a
$\K$-embedding
that extends $f_j$ as required.

One can easily verify that $\bigcup_{i<\mu}f_i$ is an embedding
of $N_2$ over $N_1$ into $M$.
\end{proof}

The following is now immediate

\begin{corollary}
Let $\K$ be an AEC without maximal models and suppose it is
Galois-stable in $\mu$.  If $\K$ has the  amalgamation property and
arbitrarily large models  then for every
$M\in\K_\mu$ and every $\sigma<\mu^+$ there exists $N\in \K_\mu$ a 
$(\mu,\sigma)$-limit model
over $M$.
\end{corollary}

\begin{remark}
In a preliminary draft of this paper we included the $\mu$-Disjoint
Amalgamation Property as one of the assumptions of Theorem \ref{exist
univ}. We thank John Baldwin and Alexei Kolesnikov for analyzing our proof
and bringing to our attention that we actually didn't use disjoint
amalgamation.

\end{remark}

%
%
%
%\begin{remark}
%For first-order complete $T$ using Robinson's consistency lemma $\Mod(T)$
%has the amalgamation property
%for all $\mu\geq |T|$.  By a theorem  of Lascar and Poizat from
%\cite{LaPo} if $T$ is
%a complete first-order theory then $\Mod(T)$ has the $\mu$-disjoint
%amalgamation property for all
%$\mu\geq |T|$.  See also Theorem 6.4.3 in Hodges \cite{Ho}.
%\end{remark}
%
%
%
%
%It is natural to ask:
%
%
%\begin{question}
%For $\K$ with the AP, is the extension
%property of Galois-types true?  More precisely can we derive
%``$\K^3_\mu$ does not have a maximal triple'' from just no maximal models
%+ AP?
%\end{question}
%
%
%%\begin{conjecture}[NEW]
%%Let $\K$ be an AEC with AP and suppose it is $\mu$-Galois-stable.  TFAE
%%\begin{enumerate}
%%
%%\item
%%$\K$ has the $\mu$-DAP.
%%
%%\item
%%Every element of $\K_\mu$ has a proper extension.
%%
%%\item
%%Every $M\in \K_\mu$ has a universal extension.
%%
%%\item
%%Dividing (splitting?) has the extension property (for $\mu$-types (see 
%Definition
%%\ref{def: extprop} below)).
%%
%%\item
%%$\K_{\mu^+}\neq\emptyset$.
%%
%%\end{enumerate}
%%\end{conjecture}
%%Notice that $(1)\implies(2)\iff(5)$ is easy
%
%

%%%%%%%%%%%%%%%%%%%%%%%%%%%%%%%%%%%%%%%%%%%%%%%%%%%%%%%%%%%%%
\bigskip
\section{Tame Abstract Elementary Classes} \label{s:tame aec}
In AECs we view types as equivalences classes of triples
(Definition \ref{defn of types}) instead of sets of formulas.  However, a
localized equivalence  relation (denoted
$E_\mu$)  is eventually utilized in various partial solutions to
Shelah's Categoricity Conjecture (see \cite{Sh 394} and \cite{Sh 576}).

 Shelah
identified
$E_\mu$  in
\cite{Sh 394}.  Here we recall the definition.

\begin{definition}
Triples $(\bar a_1, M, N_1)$ and $(\bar a_2, M, N_2)$ are said
to be \emph{$E_\mu$-related} \index{$E_\mu$-related}provided that for
every
$M'\prec_{\K}M$ with $M'\in\K_{\leq\mu}$,
$$(\bar a_1, M', N_1)E(\bar a_2, M', N_2).$$
\end{definition}

Notice that in first order logic, the finite character of consistency
implies that two types are equal if and only if they are
$E_{<\omega}$-related.

Main Claim 2.3 of part II on page 288 \cite{Sh 394} Shelah ultimately 
proves that,
under categoricity in some $\lambda>\beth_{(2^{\Hanf(\K)})^+}$ and under
the assumption that $\K$ has
the amalgamation property,
$E$-equivalence is
the same as $E_\mu$ equivalence for $1$-types over saturated models, for
some
$\mu<\lambda$.

We now define a context for abstract elementary classes
where consistency has small character.

\begin{definition}\label{tame defn}
Let $\chi$ be a cardinal number.
For $\beta>0$, we say $\K$ is \emph{$\chi$-tame for $\beta$-types}
provided that for every $M\in\K_{> \chi}$,
$p\neq  q\in \gaS^\beta(M)$  implies existence of $N\prec_{\K}M$
of cardinality $\chi$ such that  $ p\restriction N\neq q\restriction N$.
We say $\K$  is
\emph{$\chi$-tame}\index{$\chi$-tame} if it is $\chi$-tame for $1$-types.
\end{definition}

Notice that in \cite{Sh 394} 
Shelah defines a relative to tameness which considers only those types over
saturated models.  Here we are interested in types over arbitrary models.

\begin{definition}
\emph{$\K$ is tame}\index{tame}
iff there exists a $\chi<\Hanf(\K)$ such that $\K$ is $\chi$-tame.
\end{definition}

The relationship between  tameness and $\chi$-tameness is not completely
clear.  It is plausible that the following holds:

\begin{conjecture}
Let $\K$ be an AEC.  If $\K$ is $\chi$-tame for some $\chi$ then it is tame.
\end{conjecture}

\begin{remark}In this paper,
we will actually only use that $E$-equivalence is the
same as $E_\chi$-equivalence for types over limit models.
\end{remark}

Notice that if $\K$ is a finite diagram in the sense of \cite{Sh3} (i.e. we have amalgamation not
only over all models but also over all subsets of every model), then
$\K$ is a tame AEC.

  There are tame AECs with amalgamation which are not
finite diagrams.  In fact Leo Marcus in \cite{Ma} constructed
an $L_{\omega_1,\omega}$-sentence which is categorical in every
cardinality  but does not have an uncountable sequentially homogeneous
model. Recently Boris Zilber found a mathematically more natural example
\cite{Zi} motivated by Schanuel's Conjecture.  His example is not
homogeneous nor $L_{\omega_1,\omega}$-axiomatizable.

%\begin{fact}[Grossberg and Kolesnikov \cite{GrKo}]
% Let $\K$ be an AEC and $\chi\geq\LS(\K)$.
%If $\K$ has the $(\chi,3)$-amalgamation property then $\K$ is
%$\chi$-tame.
%\end{fact}
%
Shelah proves
excellence (from {\bf V}$=${\bf L}) for countable 
$L_{\omega_1,\omega}$-theories
 satisfying
$I(\aleph_{n+1},\K)<2^{\aleph_{n+1}}$ for all $n<\omega$ (\cite{Sh 87a}
and \cite{Sh 87b}). 

 In \cite{GrKo} it is shown that excellence implies
that $\K$ is $\aleph_0$-tame.  In fact, a certain two dimensional property weaker  than $(<\chi,\Pm(3))$-AP already
implies $\chi$-tameness.  Even more examples of tame classes which are not 
finite diagrams can be found from the work of Villaveces and Zambrano
(see \cite{ViZa} or \cite{Ba}).

Let $\K$ be an AEC, suppose $\kappa\geq \LS(\K)$ is a strongly compact cardinal.  If 
$\K$ is categorical in some $\lambda^+>\beth_{(2^\kappa)^+}$ then by Proposition 1.13 of \cite{Sh285}
$\K_{\geq\kappa}$ has the amalgamation property and since Galois types can be identified with sets of formulas in $L_{\kappa,\kappa}$ (see \cite{Sh285}).  Thus classes satisfying the above assumptions (categoricity above the Hanf number of a strongly compact cardinal), by Makkai and Shelah's results fit into our framework of $\kappa$-tame classes.

While we are convinced that there  are examples of arbitrary levels
of tameness at the time of submission of this paper we did not have any.

\begin{question}
For $\mu_1<\mu_2<\beth_{\omega_1}$, find an AEC which
is $\mu_2$-tame but not $\mu_1$-tame.
\end{question}

Since the submission of our paper, some progress on this question has been 
made.  Back in a discussion with Baldwin in Bogot\'{a} in 2001 we 
suggested that perhaps a non-free but almost free group of cardinality 
$\chi$ could be used to find an example of a class which is not 
$\chi$-tame.  In August 2005 Baldwin informed us that he and Shelah have a 
work in progress for class
which is not $\aleph_1$-tame they obtain this with groups and 
short exact sequences, by essentially coding the functor 
$\operatorname{Ext}(G,\mathbf Z)$ and using an $\aleph_1$-free group of 
cardinality $\aleph_1$ which is not free and not Whitehead (see 
\cite{BaSh}).

%%%%%%%%%%%%%%%%%%%%%%%%%%%%%%%%%%%%%%%%%%%%%%%%%%%%%%%%%%%%%%%%%%
\bigskip
\section{Bound on $\kappa_\mu^\beta$} \label{s:bound on kappa}

In this section we study a non-elementary relative to the first
order notion of splitting and we derive bounds for an invariant that
corresponds to $\kappa(T)$ under the assumption of stability and tameness.

First let us recall the notion of $\mu$-splitting introduced by Shelah in
\cite{Sh 394}.

\begin{definition}\index{$\mu$-splits}\index{Galois-type!$\mu$-splits} A
type
$p\in \gaS^\beta(N)$ \emph{$\mu$-splits} over $M\prec_{\K}N$ if and only if
$\|M\|\leq\mu$, there exist $N_1,N_2\in\K_{\leq\mu}$ and $h$, a
$\K$-embedding such  that $M\prec_{\K}N_l\prec_{\K}N$ for $l=1,2$ and
$h:N_1\rightarrow N_2$ such that $h\restriction M=\id_M$ and
$p\restriction N_2\neq h(p\restriction N_1)$.
\end{definition}
                        
\begin{remark}
If $T$ is a first-order theory stable in $\mu$ and $M$ is saturated, then
for all $N\prec M$ of cardinality $\mu$, the first order type,
$\ftp(a/M)$, does not split (in the first order sense) over $N$ iff
$\tp(a/M)$ does not $\mu$-split over $N$.
\end{remark}
                                                        
Notice that non-splitting is monotonic:  I.e.
If $p\in \gaS(N)$ does not $\mu$-split over $M$ (for some $M\prec_{\K}N)$
then
$p$ does not $\mu$-split over $M'$ for every $M\prec_{\K}M'\prec_{\K}N$ of
cardinality $\mu$.

Similar to $\kappa(T)$ when $T$ is first-order the following is a
natural cardinal
invariant of $\K$:

\begin{definition}
Let $\beta>0$.
We define $\kappa^\beta_\mu(\K)$
\index{$\kappa^\beta_\mu(\K)$}to be the minimal $\kappa<\mu^+$ such that
for every
$\langle M_i\in \K_{\mu}
\mid i\leq\kappa\rangle$ which satisfies
\begin{enumerate}
\item $\kappa=\cf(\kappa)<\mu^+$,
\item $\langle M_i\mid i\leq\kappa\rangle$ is $\prec_{\K}$-increasing and
continuous and
\item for every $i<\kappa$, $M_{i+1}$ is a $(\mu,\theta)$-limit over
$M_i$ for some $\theta<\mu^+$,
\end{enumerate}
and for every $p\in \gaS^\beta(M_\kappa)$,
there exists $i<\kappa$ such that $p$ does not $\mu$-split over $M_i$.
If no such $\kappa$ exists, we say
$\kappa^\beta_{\mu}(\K)=\infinity$.
\end{definition}

Another relative of $\kappa(T)$ is the following
\begin{definition}\index{$\bar \kappa^\beta_\mu(\K)$}
For $\beta>0$, $\bar \kappa^\beta_\mu(\K)$ is the minimal cardinal $\bar
\kappa$ such that
for every $N\in\K_\mu$ and every $p\in\gaS^\beta(N)$ there are
$\lambda<\bar
\kappa$ and $M\in\K_\lambda$ such that $p$ does not $\lambda$-split over
$M$.
\end{definition}

It is not difficult to verify that
\begin{proposition}For $\mu$ with $\cf(\mu)>\bar\kappa^\beta_\mu(\K) $,
we have
$\kappa^\beta_\mu(\K)\leq\bar\kappa^\beta_\mu(\K)$.
\end{proposition}

While for stable first order theories (when $\beta<\omega$) both
invariants are equal, the situation for non-elementary classes is more
complicated.  Already in \cite{Sh 394} it was observed that
$\kappa^\beta_\mu(\K)$ is better behaved than $\bar\kappa^\beta_\mu(\K)$
when a bound for $\kappa^\beta_\mu(\K)$ was found.  A corresponding bound
for $\bar\kappa^\beta_\mu(\K)$ is unknown.  We will defer dealing with
the  invariant $\bar\kappa^\beta_\mu(\K)$ to a
future paper.

Notice that Theorem 2.2.1 of \cite{ShVi 635}
states that certain categorical abstract elementary classes
satisfy
$\kappa^1_\mu(\K)=\omega$, for various $\mu$.
A slight modification of the argument of Claim 3.3 from \cite{Sh 394} can
be used to prove a related result using the weaker assumption of
Galois-stability only:

\begin{fact}\label{bounded kapp}
Let $\beta>0$.  Suppose that $\K$ is
$\beta$-Galois-stable in $\mu$.
For every $p\in \gaS^\beta(N)$ there exists
$M\prec_{\K}N$ of cardinality $\mu$
such that $p$ does not $\mu$-split over $M.$
Thus $\bar\kappa^\beta_\mu(\K)\leq\mu$.
\end{fact}

For the sake of completeness an argument for Fact \ref{bounded kapp}
is included:
\begin{proof}
Suppose  $p=\tp(\bar a/N,\C)$ with $\bar a\in\sq\beta|\C|$
and
$p$ $\mu$-splits over $M$, for every $M\prec_{\K}N$ of cardinality
$\mu$.  We may assume that $\bar a\notin N$, since it is routine to
check that an algebraic type never splits over its domain.

Let $\chi:=\min\{\chi\mid 2^\chi>\mu\}$.  Notice that $\chi\leq
\mu$ and
$2^{<\chi}\leq\mu$.

First we will define $\langle M_\alpha,N_{\alpha,\ell}\prec
N\mid\alpha<\chi,\;\ell=1,2\rangle\subseteq
\K_\mu$  which will then be used to construct
another sequence $\langle M^*_\alpha\in\K_\mu\mid\alpha\leq\chi\rangle$ so
that
\[
|\gaS^\beta(M^*_\chi)|\geq
2^\chi>\mu \text{ will witness
 }
\text{Galois-instability in }\mu.
\]

The definition of $\langle M_\alpha,N_{\alpha,1},N_{\alpha,2}\mid
\alpha<\chi\rangle$ follows by induction on $\alpha$ and our choice of
$p$.  Let
$M_0\prec M$ be any model of cardinality
$\mu$.

For $\alpha=\gamma+1$, by our choice of $p$, we know $p$ must $\mu$-split
over $M_\gamma$.
Then, there are
$N_{\gamma,\ell}\prec_{\K} M $ of cardinality $\mu$ for $\ell=1,2$
and there is
$F_\gamma:N_{\gamma,1}\cong_{M_\gamma}N_{\gamma,2}$ such that
$F_\gamma(p\restriction N_{\gamma,1})\neq p\restriction
N_{\gamma,2}$.  Pick $M_\alpha\prec_{\K}M$ of cardinality $\mu$ containing
the set
$|N_{\gamma,1}|\cup |N_{\gamma,2}|$.

When $\alpha$ is a limit ordinal we define
$M_\alpha:=\Union_{\gamma<\alpha}M_\gamma$.

 Now
we use our sequence $\langle M_\alpha\mid\alpha<\chi\rangle $ to define
for $\alpha\leq\chi$ another sequence of models
$M_\alpha^*\in \K_\mu$ and a tree of $\K$-embeddings
$\langle h_\eta\mid\eta\in
\sq{\alpha}2\rangle$ satisfying
\begin{enumerate}

\item
$\gamma<\alpha\implies M^*_\gamma\prec_{\K} M^*_\alpha$,

\item
for $\alpha$ limit let $M^*_\alpha=\bigcup_{\gamma<\alpha}M^*_\gamma$,

\item
$\gamma<\alpha\wedge \eta\in \sq{\alpha}2\implies h_{\eta\restriction
\gamma}\subseteq h_\eta$,

\item
$\eta\in\sq{\alpha}2\implies h_\eta:M_\alpha \rightarrow
M^*_\alpha$
and

\item
$\alpha=\gamma+1\wedge \eta\in \sq{\alpha}2\implies
h_{\eta\conc 0}(N_{\gamma,1})=h_{\eta\conc 1}(N_{\gamma,2})$

\end{enumerate}

The construction is possible by using the $\mu$-amalgamation
property  several times.
For $\alpha=0$, define $M^*_0=M$ and
$h_{\langle\rangle}:=\id_{M^*_{0}}$.  When $\alpha$ is a limit
ordinal we define $M^*_\alpha:=\Union_{\gamma<\alpha}M_\gamma$ and for
$\eta\in\sq\alpha 2$,
$h_\eta:=\Union_{\gamma<\alpha}h_{\eta\restriction\gamma}$.

Suppose we have completed the construction for $\gamma$ and
$\alpha=\gamma+1$.   Given $\eta\in \sq{\gamma}2$, let $\bar h_\eta$ be
an extension of $h_\eta$ to an automorphism of $\C$.  Notice that $\bar
h_\eta\circ F_\gamma(N_{\gamma,1})=\bar h_{\eta}(N_{\gamma,2})$ by our
choice of $F_\gamma$, $N_{\gamma,1}$ and $N_{\gamma,2}$.  Let
$h_{\eta\conc\langle 0\rangle}$ be some extension of
$(\bar h_\eta\restriction N_{\gamma,2})\circ F_{\gamma}$ to $M_{\gamma+1}$
and let
$h_{\eta\conc\langle 1\rangle}$ be $\bar h_{\eta}\restriction
M_{\gamma+1}$.  Let $M^*_{\gamma+1}$ be a model of cardinality $\mu$
extending $M^*_\gamma$ and containing $h_{\eta\conc\langle
0\rangle}(M_{\gamma+1})$ and $h_{\eta\conc\langle
1\rangle}(M_{\gamma+1})$ for each $\eta\in\sq\gamma 2$.

This completes the construction and for each $\eta\in\sq\chi 2$ we can
fix $H_\eta\in\Aut(\C)$ extending $h_\eta$.

\begin{claim}\label{many types claim}
$\eta\neq \nu\in \sq{\chi}2\implies
\tp(H_\eta(\bar a)/M^*_\chi)\neq
\tp(H_\nu(\bar a)/M^*_\chi)$.
\end{claim}
\begin{proof}
Let $\rho:=\eta\wedge\nu$ and suppose that $\rho\conc \langle
0\rangle<\eta$ and
$\rho\conc\langle 1\rangle<\nu$.
Let $\gamma$ be the length of $\rho$.
Suppose for the sake of contradiction that $$
\tp(H_\eta(\bar a)/M^*_\chi)=\tp(H_\nu(\bar a)/M^*_\chi).$$

In particular,
$$\tp(H_\eta(\bar a)/H_{\rho\conc\langle
0\rangle}(N_{\gamma,1}))=\tp(H_\nu(\bar a)/H_{\rho\conc\langle
1\rangle}(N_{\gamma,2})).$$

Referring back to our choice of $\bar h_\rho$ and $h_{\rho\conc \langle
0\rangle}$, we can apply $H_{\nu}^{-1}$ to obtain

$$\tp(H_\nu^{-1}\circ H_\eta(\bar a)/F_\gamma(N_{\gamma,1})=
\tp(\bar a/N_{\gamma,2}).$$

All we need to verify in order to contradict our original choices of
$N_{\gamma,1}$, $N_{\gamma,2}$ and $F_\gamma$ is that
for some extension of $F_\gamma$ to an automorphism of $\C$ we have that
$$(*)\quad\tp(F_\gamma(\bar a)/N_{\gamma,2})=\tp(H_\nu^{-1}\circ
H_\eta(\bar a)/N_{\gamma,2}).$$

Notice that by our choice of $H_\nu$ and $H_\eta$ we have that
$H_\nu^{-1}\circ H_\eta\restriction N_{\gamma,1}=F_\gamma$.  Thus
$H_\nu^{-1}\circ H_\eta\supseteq F_\gamma$, so by definition of the image
of a type,
$F_\gamma(p\restriction N_{\gamma,1})=\tp(H_\nu^{-1}\circ H_\eta(\bar
a)/N_{\gamma,2})$.  Combining this with $(*)$ we obtain a contradiction
to our choice of $N_{\gamma,1}$, $N_{\gamma,2}$ and $F_{\gamma}$
which witness the splitting of $p$ over $M_\gamma$.
\end{proof}
Thus Claim \ref{many types claim} gives us that $|\gaS(M^*_\chi)|\geq
2^\chi>\mu$ and
$\K$ is not
$\beta$-Galois-stable in
$\mu$.
\end{proof}

In Theorem \ref{tame and stable imply kappa exists} below we present an
improvement of Fact  \ref{bounded kapp}
for tame AECs:  In case $\K$ is $\beta$-stable in $\mu$
for some $\mu$ above its second Hanf number then
$\kappa^\beta_\mu(\K)$ is bounded by the  second Hanf number.  Notice
that the bound does not depend on $\mu$.

We will need to relate Galois-stability with the failure of a certain
infinitary order property in order to find a bound for
$\kappa^\beta_{\mu}(\K)$. The order property, defined next, is an analog
of the first order definition using formulas.  This order property for
non-elementary classes  was introduced by Shelah in
\cite{Sh 394}.

\begin{definition}\label{op}
$\K$ is said to have the
\emph{$\chi$-order property}\index{$\kappa$-order property} provided
that for every $\alpha$, there exists $\langle \bar d_i\mid
i<\alpha\rangle$ and  where
$\bar d_i\in\sq\chi\C$ such that if
$i_0<j_0<\alpha$ and $i_1<j_1<\alpha$,
$$\quad(*)\text{ then for
no }f\in Aut(\C)\text{ do we have }f(\bar d_{i_0}\conc\bar d_{j_0})
=\bar d_{j_1}\conc\bar d_{i_1}.$$
\end{definition}

\begin{example}
Let $T$ be a  first order theory in a countable language.
Set $\K:=\Mod(T)$ and $\prec_{\K}$ the usual elementary submodel
relation.  $\K$ has the $\aleph_0$-order property iff $T$ is unstable of 
$T$ has dop
or otop.
\end{example}

\begin{remark}[Trivial monotonicity]
Notice that for $\kappa_1<\kappa_2$ if a class has the $\kappa_1$-order
property then it has the $\kappa_2$-order property.
\end{remark}

\begin{fact}[Claim 4.6.3 of \cite{Sh 394}]\label{Claim 4.6.3}
We may replace the phrase \emph{every $\alpha$} in Definition \ref{op} with
\emph{every $\alpha <\beth_{(2^{\chi+\LS(\K)})^+}$} and get an
equivalent definition.
\end{fact}

\begin{fact}[Claim 4.8.2 of \cite{Sh 394}]\label{op implies instability}
If $\K$ has the $\chi$-order property and $\mu\geq\chi$, then for
some $M\in\K_\mu$ we have that $|\gaS^\chi(M)/E_\chi|\geq\mu^+$.
Moreover, we can conclude that $\K$ is not $\chi$-Galois-stable in $\mu$.
\end{fact}

The following is a generalization of an old theorem (Theorem 4.17 in
\cite{GrLe2}) of Shelah from
 \cite{Sh3}.  A special case of this theorem is Theorem 
\ref{bounded non-split} in the
abstract.
\begin{theorem}\label{tame and stable imply kappa exists}
Let $\beta>0$.
Suppose that $\K$ is $\chi$-tame for $\beta$-types.   If $\K$ is
$\beta$-stable in $\mu$ with
$\beth_{(2^{\chi+\LS(\K)})^+}\leq\mu$, then
$\kappa_\mu^\beta(\K)<\beth_{(2^{\chi+\LS(\K)})^+}$.
\end{theorem}
\begin{proof}
Let $\chi':=\beth_{(2^{\chi+\LS(\K)})^+}$.
Suppose that the conclusion of the theorem does not hold.
Let
$\langle M_i\in \K_{\mu}
\mid i\leq\chi'\rangle$ and $p\in \gaS^\beta(M_{\chi'})$
witness the failure.  Namely, the following hold:
\begin{enumerate}
\item $\langle M_i\mid i\leq\chi'\rangle$ is $\prec_{\K}$-increasing and
continuous,
\item for every $i<\chi'$, $M_{i+1}$ is a $(\mu,\theta)$-limit over
$M_i$ for some $\theta<\mu^+$ and
\item\label{last item}
for every $i<\mu^+$, $p$ $\mu$-splits over $M_i$.
\end{enumerate}

For every $i<\chi'$ let $f_i, N_i^1$ and $N_i^2$ witness that $p$
$\mu$-splits over $M_i$.  Namely,

$$M_i\prec_{\K}N_i^1,N_i^2\prec_{\K}M,$$
$$f_i:N_i^1\cong N_i^2\text{ with }f_i\restriction M_i=\id_{M_i}$$
$$\text{ and }f_i(p\restriction N_i^1)\neq p\restriction N_i^2.$$

By $\chi$-tameness, there exist models $B_i$ and $A_i:=f^{-1}_i(B_i)$ of
size
$<\chi$ such that
$$f_i(p\restriction A_i)\neq p\restriction B_i.$$

By renumbering our chain of models, we may assume that
\begin{enumerate}\setcounter{enumi}{3}
\item $A_i,B_i\subset M_{i+1}$.
\end{enumerate}

Since $M_{i+1}$ is a limit model over $M_i$, we can additionally fix
\begin{enumerate}\setcounter{enumi}{4}
\item $\bar c_i\in M_{i+1}$ realizing $p\restriction M_i$.
\end{enumerate}

For each $i<\mu$, let $\bar d_i:=A_{i}\conc B_{i}\conc \bar c_{i}$.

\begin{claim}
$\langle \bar d_i\mid i<\chi'\rangle$ witnesses the $\chi$-order
property.
\end{claim}
\begin{proof}
Suppose for the sake of contradiction that there exist $g\in Aut(\C)$,
$i_0<j_0<\chi'$ and $i_1<j_1<\chi'$ such that
$$g(\bar d_{i_0}\conc\bar d_{j_0})
=\bar d_{j_1}\conc\bar d_{i_1}.$$

Notice that since $i_0<j_0<\alpha$ we have that $\bar c_{i_0}\in
M_{j_0}$. So $f_{j_0}(\bar c_{i_0})=\bar c_{i_0}$.  Recall that
$f_{j_0}(A_{j_0}) =B_{j_0}$.
Thus, $f_{j_0}$ witnesses that
$$(*)\quad\tp(\bar c_{i_0}\conc A_{j_0}/\emptyset)
=\tp(\bar c_{i_0}\conc B_{j_0}/\emptyset).$$

Applying $g$ to $(*)$ we get
$$(**)\quad\tp(\bar c_{j_1}\conc A_{i_1}/\emptyset)
=\tp(\bar c_{j_1}\conc B_{i_1}/\emptyset).$$

Applying $f_{i_1}$ to the RHS of $(**)$, we notice that
$$(\sharp)\quad\tp(f_{i_1}(\bar c_{j_1})\conc B_{i_1}/\emptyset)
=\tp(\bar c_{j_1}\conc B_{i_1}/\emptyset).$$

Because $i_1<j_1$, we have that
$\bar c_{j_1}$ realizes $p\restriction M_{i_1}$.
Thus, $(\sharp)$ implies
$$f_{i_1}(p\restriction A_{i_1})
=p\restriction B_{i_1},$$
which contradicts our choice of $f_{i_1}$, $A_{i_1}$ and $B_{i_1}$.

\end{proof}

By Claim \ref{Claim 4.6.3} and Fact \ref{op implies instability}, we
have that $\K$ is $\beta$-unstable in
$\mu$, contradicting our hypothesis.

\end{proof}

%%%%%%%%%%%%%%%%%%%%%%%%%%%%%%%%%%%%%%%%%%%%%%%%%%%%%%%%%%%%%%%%%%
\bigskip
\section{Morley sequences} \label{s:main theorem}

We have now developed enough splitting machinery to
derive the existence of Morley sequences for tame, Galois-stable AECs.

The following is a new Galois-theoretic notion of indiscernible sequence.

\begin{definition}
\begin{enumerate}
\item
$\langle \bar a_i\mid i<i^*\rangle$ is a  \emph{Galois-indiscernible
sequence over
$M$} iff for every $i_1<\dots<i_n<i^*$ and every $j_1<\dots<j_n<i^*$,
$\tp(\bar a_{i_1}\dots\bar a_{i_n}/M)=\tp(\bar a_{j_1}\dots\bar
a_{j_n}/M)$.
\item
$\langle \bar a_i\mid i<i^*\rangle$ is a \emph{Galois-indiscernible
sequence over
$A$}\index{Galois-indiscernible sequence} iff for every
$i_1<\dots<i_n<i^*$ and every
$j_1<\dots<j_n<i^*$, there exists $M_i, M_j, M^*\in \K$ and
$\prec_{\K}$-mappings $f_i,f_j$  such that
\begin{enumerate}
\item
$A\subseteq M_i, M_j$;
\item
$f_\ell:M_\ell\rightarrow M^*$, for $\ell\in \{i,j\}$;
\item
$f_i(\bar a_{i_0},\dots,\bar a_{i_n})=
f_j(\bar a_{j_0},\dots,\bar a_{j_n})$ and
\item and $f_i\restriction A=f_j\restriction A = \id_A$.
\end{enumerate}
\end{enumerate}
\end{definition}

\begin{remark}
This is on the surface a weaker notion of indiscernible sequence than is
presented in
\cite{Sh 394}.  However, under
the amalgamation property, this definition and
the definition in \cite{Sh 394} are equivalent.
\end{remark}

Lemma \ref{ind} provides us with sufficient conditions to find
an indiscernible sequence.

\begin{lemma}\label{ind}
Let $\mu\geq \LS(\K)$, $\kappa,\lambda$ be ordinals and
$\beta$ a positive ordinal.
Suppose that $\langle M_i\mid i<\lambda\rangle$ and
$\langle \bar a_i\mid i<\lambda\rangle$ satisfy
\begin{enumerate}
\item
$\langle M_i\in\K_{\mu}\mid i<\lambda\rangle$ are
$\prec_{\K}$-increasing;

\item\label{universal}
$M_{i+1}$ is a $(\mu,\kappa)$-limit over $M_i$;
\item
$\bar a_i\in \sq{\beta}M_{i+1}$;
\item\label{non-split}
$p_i:=\tp(\bar a_i/M_i,M_{i+1})$ does not $\mu$-split over $M_0$ and
\item\label{inc}
for $i<j<\lambda$, $p_i\leq p_j$.
\end{enumerate}
Then,
$\langle \bar a_i\mid i< \lambda\rangle$ is a Galois-indiscernible
sequence over $M_ 0$.
\end{lemma}

\begin{definition}
A sequence $\langle \bar a_i,M_i\mid i<\lambda\rangle$ satisfying conditions
$(1)-(\ref{inc})$ of Lemma \ref{ind} is called a
\emph{Morley sequence}.  It is a \emph{Morley sequence for 
$p\in\gaS(\Union_{i<\lambda}M_i)$} if $p$ extends each $\tp(\bar a_i/M_i)$.
For $A\subseteq M_0$, such a sequence  is a \emph{Morley sequence over
$A$}.
\end{definition}

\begin{remark}
Notice that our definition of Morley sequence varies from some
literature.  An alternative name for our sequences (suggested by John
Baldwin) is a \emph{coherent non-splitting sequence}.\index{coherent
non-splitting sequence}
\end{remark}

\begin{remark}
While the statement of the lemma is similar to Shelah's first order Lemma
I.2.5 in \cite{Shc}, the proof differs, since types are not
sets of formulas.
\end{remark}

\begin{proof}
We prove that for $i_0<\dots<i_n<\lambda$ and $j_0<\dots<j_n<\lambda$,
$\tp(\bar a_{i_0},\dots,\bar a_{i_n}/M_0,M_{i_{n+1}})=
\tp(\bar a_{j_0},\dots,\bar a_{j_n}/M_0,M_{j_{n+1}})$ by induction
on $n<\omega$.

$n=0$:  Let $i_0,j_0<\lambda$ be given.  Condition \ref{inc}, gives us
$$\tp(\bar a_{i_0}/M_0,M_{i_0+1})=
\tp(\bar a_{j_0}/M_0,M_{j_0+1}).$$

$n>0$: Suppose that the claim holds for all increasing
sequences $\bar i$ and $\bar j\in \lambda$ of length $n$.  Let
$i_0<\dots<i_{n}<\lambda$ and $j_0<\dots<j_{n}<\lambda$ be given.
Without loss of generality, $i_n\leq j_n$.
Define $M^*:=M_1$.
 From condition \ref{universal} and uniqueness of $(\mu,\omega)$-limits,
we can find
a $\prec_{\K}$-isomorphism, $g:M_{j_n}\rightarrow
M_{i_n}$ such that $g\restriction M_{0}=\id_{M_{0}}$.
Moreover we can extend $g$ to $g:M_{j_n+1}\rightarrow M_{i_n+1}.$
Denote by $\bar b_{j_l}:=g(\bar a_{j_l})$ for $l=0,\dots,n$.  Notice that
$b_{j_l}\in M_{i_n}$ for $l<n$.  Since
$\tp(\bar b_{j_0},\dots,\bar b_{j_n}/M_0,M_{i_n+1})=\tp(\bar
a_{j_0},\dots, \bar a_{j_n}/M_0,M_{j_n+1})$
it suffices to prove that
$\tp(\bar b_{j_0},\dots,\bar b_{j_n}/M_0,M_{i_n+1})=
\tp(\bar a_{i_0},\dots,\bar a_{i_n}/M_0,M_{i_n+1}).$

Also notice that our $\prec_{\K}$-mapping $g$ preserves some properties
of
$p_j$.  Namely, since $p_j$ does not $\mu$-split over $M_0$,
$g(p_j\restriction M_{j_n})=p_j\restriction M_{i_n}$. \\
 Thus,
$\tp(\bar b_{j_n}/M_{i_n},M_{i_n+1})=
\tp(\bar a_{j_n}/M_{i_n},M_{i_n+1}).$
In particular we have that
$\tp(\bar b_{j_n}/M_{i_n},M_{i_n+1})$ does not
$\mu$-split over $M_0$.

By the induction
hypothesis
$$\tp(\bar b_{j_0},\dots,\bar b_{j_{n-1}}/M_0,M_{i_n})=\tp(\bar
a_{i_0},\dots, \bar a_{i_{n-1}}/M_0,M_{i_n}).$$
Thus we can find $h_i:M_{i_n+1}\rightarrow M^*$ and
$h_j:M_{i_n+1}\rightarrow M^*$ such that
$h_i(\bar a_{i_0},\dots,\bar a_{i_{n-1}})=
h_j(\bar b_{j_0},\dots,\bar b_{j_{n-1}})$.  Let us abbreviate
$\bar b_{j_0},\dots,\bar b_{j_{n-1}}$ by $\bar b_{\bar j}$.
Similarly we will write $\bar a_{\bar i}$ for
$\bar a_{i_0},\dots,\bar a_{i_{n-1}}$.

By appealing to condition \ref{non-split}, we derive several equalities that
will be useful in the latter portion of the proof.
Since
$p_j$ does not
$\mu$-split over $M_0$, we have that
$p_j\restriction h_j(M_{i_n})=h_j(p_j\restriction M_{i_n})$, rewritten as
$$(*)\quad \tp(\bar b_{j_n}/h_j(M_{i_n}),M_{i_n+1})=
\tp(h_j(\bar b_{j_n})/h_j(M_{i_n}),M^*).$$

Similarly as $p_i$ does not $\mu$-split over $M_0$, we get\\
$p_i\restriction h_j(M_{i_n})=h_j(p_i\restriction M_{i_n})$
and $p_i\restriction h_i(M_{i_n})=h_i(p_i\restriction M_{i_n})$.
These equalities translate to
$$(**)_j\quad \tp(\bar a_{i_n}/h_j(M_{i_n}),M_{i_n+1})=
\tp(h_j(\bar a_{i_n})/h_j(M_{i_n}),M^*)\text{ and}$$
$$(**)_i\quad \tp(\bar a_{i_n}/h_i(M_{i_n}),M_{i_n+1})=
\tp(h_i(\bar a_{i_n})/h_i(M_{i_n}),M^*),\text{ respectively}.$$

Finally, from condition \ref{inc}, notice that
$$(***)\quad \tp(\bar a_{i_n}/M_{i_n},M_{i_n+1})=
\tp(\bar b_{j_n}/M_{i_n},M_{i_n+1}).$$

%Let $N^*_1:=h_i(M_{i_n})$ and
%$N^*_2:=h_j(M_{i_n})$.  Consider
%$H^*:=h_j\circ h_i^{-1}:N^*_1\rightarrow N^*_2$ and
%the type $p^*:=p_{j_n}\restriction M_{i_n}=\tp(\bar
%b_{j_n},M_{i_n},M_{i_n+1})$.  Since $p_{j_n}$ does not
%$\mu$-split over $M_0$, we have that
%$p^*\restriction N_2^*=H(p^*\restriction N_1^*)$.
%
%First notice that
%$p^*\restriction N_2^*=\tp(\bar b_{j_n},h_j(M_{i_n}),M_{i_n+1})$.
%Moreover $(*)$ gives us that
%$$(\sharp)\quad  p^*\restriction N_2^*=\tp(h_j(\bar b_{j_n}),
%h_j(M_{i_n}),M^*).$$
%
%Now consider $H(p^*\restriction N_1^*)$.
%$p^*\restriction N_1^*=\tp(\bar b_{j_n},h_i(M_{i_n}),M_{i_n+1})$.
%Since $h_i(M_{i_n})\prec_{\K}M_{i_n}$, equality $(***)$ yields
%$p^*\restriction N_1^*=\tp(\bar a_{i_n},h_i(M_{i_n}),M_{i_n+1})$.
%Furthermore, we can apply equality $(**)_i$ to conclude
%$p^*\restriction N_1^*=\tp(h_i(\bar a_{i_n}),h_i(M_{i_n}),M^*)$.
%Thus an application of $H$ yields
%$$(\sharp\sharp)\quad
%H(p^*\restriction N_1^*)=
%\tp(h_j(\bar a_{i_n}),h_j(M_{i_n}),M^*).$$
%
%Therefore combining equalities $(\sharp)$ and $(\sharp\sharp)$
%we get that
%$$(\dagger)\quad\tp(h_j(\bar b_{j_n}),h_j(M_{i_n}),M^*)=
%\tp(h_j(\bar a_{i_n}),h_j(M_{i_n}),M^*).$$

Applying $h_j$ to $(***)$ yields
$$(\dagger)\quad\tp(h_j(\bar b_{j_n})/h_j(M_{i_n}),M^*)=
\tp(h_j(\bar a_{i_n})/h_j(M_{i_n}),M^*).$$

Since $h_i(\bar a_{\bar i})=h_j(\bar b_{\bar j})\in h_j(M_{i_n})$,
we can draw from $(\dagger)$ the following:
$$(1)\quad\tp(h_j(\bar b_{j_n})\conc h_j(\bar b_{\bar j})/M_0,M^*)=
\tp(h_j(\bar a_{j_n})\conc h_i(\bar a_{\bar i})
/M_0,M^*).$$

Similarly, equality $(**)_i$ allows us to see
$$(2)\quad
\tp(\bar a_{i_n}\conc h_i(\bar a_{\bar i})/M_0,M^*)=
\tp(
h_i(\bar a_{i_n})\conc h_i(\bar a_{\bar i})/M_0,M^*)
.$$

Since
$\tp(h_j(\bar a_{i_n})/h_j(M_{i_n}),M^*)=
\tp(\bar a_{i_n}/h_j(M_{i_n}),M_{i_n+1})$ (equality $(**)_j)$) and
$h_i(\bar a_{\bar i})=
h_j(\bar b_{\bar j})\in
h_j(M_{i_n})$, we get that
$$(3)\quad\tp( h_j(\bar a_{i_n})\conc h_i(\bar a_{\bar i})/M_0,M^*)=
\tp(
\bar a_{i_n}\conc h_i(\bar a_{\bar i})/M_0,M^*).$$

Combining equalities $(1)$, $(2)$ and
$(3)$, we get
$$(\dagger\dagger)\quad
\tp(h_i(\bar a_{\bar i})\conc h_i(\bar a_{i_n})/
M_0,M^*)=
\tp(h_j(\bar b_{\bar j})\conc h_j(\bar b_{j_n})/
M_0,M^*).$$

Recall that $h_i\restriction M_0=h_j\restriction M_0=\id_{M_0}$.
Thus $(\dagger\dagger)$, witnesses that
$$\tp(\bar a_{i_0},\dots,\bar a_{i_n}/M_0,M_{i_n+1})=
\tp(\bar b_{j_0},\dots,\bar b_{j_n}/M_0,M_{i_n+1}).$$
\end{proof}

\begin{theorem}\label{ind thm}\index{Morley sequence!existence|(}
\index{Galois-indiscernible sequence!existence|(}
Fix $\beta>0$.  Suppose $\K$ is $\chi$-tame for $\beta$-types with
$\chi<\beth_{(2^{\Hanf(\K)})^+}$.  Suppose $\mu\geq
\beth_{(2^{Hanf(\K)})^+}$. Let $M\in\K_{>\mu}$, $A\subset M$ and
$I\subseteq\sq\beta M$ be given such that
$|I|\geq\mu^+>|A|$.
If $\K$ is Galois $\beta$-stable in $\mu$,
then there exists $J\subset I$ of cardinality $\mu^+$,
Galois indiscernible over $A$.  Moreover $J$ can be chosen to be a Morley
sequence  over $A$.
\end{theorem}

\begin{proof}
Fix $\kappa=\kappa^\beta_\mu$.  By Theorem \ref{tame and stable imply
kappa exists},
$\kappa<\beth_{(2^{\Hanf(\K)})^+}$. Let $\{\bar a_i\in I\mid i<\mu^+\}$
be given.  Define
$\langle M_i\in K_\mu\mid i<\mu^+\rangle$, $\prec_{\K}$-increasing and
continuous, satisfying
\begin{enumerate}
\item $A\subseteq |M_0|$
\item $M_{i+1}$ is a $(\mu,\kappa)$-limit over $M_i$
\item $\bar a_i\in M_{i+1}$
\end{enumerate}
Let $p_i:=\tp(\bar a_i/M_i,M_{i+1})$ for every $i<\mu^+$.  Define
$f:S^{\mu^+}_{\kappa}\rightarrow \mu^+$ by
$$f(i):=\min\{j<\mu^+\mid p_i\text{ does not }\mu\text{-split over
}M_j\}.$$ By our choice of $\kappa$,
$f$ is regressive.  Thus by Fodor's Lemma,
there are a stationary set $S\subseteq S^{\mu^+}_\kappa$ and $j_0<\mu^+$
such that for every $i\in S$,
$$(\dagger)\quad p_i \text{ does not }\mu\text{-split over }M_{j_0}.$$
By stability and the pigeon-hole
principle there exists $p^*\in \gaS(M_{j_0})$ and
$S^*\subseteq S$ of cardinality $\mu^+$ such that
for every $i\in S^*$, $p^*=p_i\restriction M_{j_0}$.
Let $M^*:=M_{j_0}$.
Enumerate
and $\langle M_{j_i}\mid i\in S^*\rangle$.
Let $M^{**}:=M_{j_1}$ in this  enumeration.  Again, by stability we can
find
$S^{**}\subset S^*$ of cardinality $\mu^+$ such that
for every $i\in S^{**}$, $p^{**}=p_i\restriction M^*$.
Notice that $M^*,M^{**}\prec_{\K}M_i$ for every $i\in S^{**}$.

\begin{subclaim}\label{inc claim}
For $i<j\in S^{**}$, $p_i=p_j\restriction M_i$.
\end{subclaim}

\begin{proof}
Let $0<i<j\in S^{**}$ be given.
Since $M_{i+1}$ and $M_{j+1}$ are $(\mu,\kappa)$-limits over
$M_i$, there exists an isomorphism $g:M_{j+1}\rightarrow M_{i+1}$
such that $g\restriction M_i=\id_{M_i}$.  Let $\bar b_j:=g(\bar a_j)$.
Since the type $p_j$ does not $\mu$-split over $M^*$, $g$
cannot witness the splitting.  Therefore, it must be the case that
$\tp(\bar b_j/M_i,M_{i+1})=p_j\restriction M_i$.
Then, it suffices to show that $p_i=\tp(\bar b_j/M_i,M_{i+1})$.

Since $p_i\restriction M^{*}= p_j\restriction M^{*}$, we can find
$\prec_{\K}$-mappings witnessing the equality.  Furthermore since
$M^{**}$ is universal over $M^*$, we can find
$h_l:M_{l+1}\rightarrow M^{**}$ such that
$h_l\restriction M^*=\id_{M^*}$ for $l=i,j$ and
$h_i(\bar a_i)=h_j(\bar b_j)$.

We will use $(\dagger)$ to derive several inequalities.  Consider the 
following
possible witness to splitting.  Let $N_1:=M_i$ and $N_2:=h_i(M_i)$.  Since 
$p_i$ does not
$\mu$-split over $M^*$, we have that $p_i\restriction N_2=
h_i(p_i\restriction N_1)$, rewritten as
$$(*)\quad \tp(\bar a_i/h_i(M_i),M_{i+1})=
\tp(h_i(\bar a_i)/h_i(M_i), M^{**}).$$

Similarly we can conclude that
$$(**)\quad \tp(\bar b_j/h_j(M_i), M_{i+1})=
\tp(h_j(\bar b_j)/h_j(M_i), M^{**}).$$

By choice of $S^{**}$, we know that
$$(***)\quad \tp(\bar b_j/M^{**})=\tp(\bar a_i/M^{**}).$$

Now let us consider another potential witness of splitting.
$N^*_1:=h_i(M_i)$ and $N^*_2:=h_j(M_i)$ with
$H^*:=h_j\circ h_i^{-1}:N^*_1\rightarrow N^*_2$.  Since
$p_j\restriction M_i$ does not
$\mu$-split over $M_0$, $p_j\restriction N^*_2 =
H^*(p_j\restriction N^*_1)$.  Thus by $(**)$ we have
$$(\sharp)\quad H^*(p_j\restriction N_1^*) = \tp(h_j(\bar
b_j)/h_j(M_i),M^{**}).$$

Now let us translate $H^*(p_j\restriction N^*_1)$.
By monotonicity and $(***)$, we have that
$p_j\restriction N^*_1=\tp(\bar b_j/h_i(M_i),M_{i+1})=
\tp(\bar a_i/h_i(M_i),M_{i+1})$.
We can then conclude by $(*)$ that
$p_j\restriction N^*_1=\tp(h_i(\bar a_i)/h_i(M_i),M_{i+1})$.
Applying $H^*$ to this equality yields
$$(\sharp\sharp)\quad H^*(p_j\restriction N^*_1)=\tp(h_j(\bar
a_i)/h_j(M_i),M^{**)}.$$

By combining the equalities from $(\sharp)$ and $(\sharp\sharp)$
and applying $h_j^{-1}$ we get that
$$\tp(\bar b_j/M_i,M_{i+1})=\tp(\bar a_i/M_i,M_{i+1}).$$
\end{proof}

Notice that by Subclaim \ref{inc claim} and our choice of $S^{**}$,
$\langle M_i\mid i\in S^{**}\rangle$ and $J:=\langle \bar a_i\mid i\in
S^{**}\rangle$  satisfy the conditions of Lemma \ref{ind}.
Applying Lemma \ref{ind}, we get that
$\langle \bar a_i\mid i\in S^{**}\rangle$ is a Morley sequence over
$M^*$.  In particular, since $A\subset M^*$, we have that
$\langle \bar a_i\mid i\in S^{**}\rangle$ is a Morley sequence over
$A$.

\end{proof}\index{Morley sequence!existence|)}
\index{Galois-indiscernible sequence!existence|)}

\begin{remark}
Our effort to prove the results of this section in summer 2001 lead us to 
the notion of tameness as defined in this paper.  A preliminary version of 
this paper containing the results of this section was posted by us on the 
web already in July 27th, 2001 with the title \emph{Morley Sequences in 
AECs}.
\end{remark}

%%%%%%%%%%%%%%%%%%%%%%%%%%%%%%%%%%%%%%%%%%%%%%%%%%%%%%%%%%%%%%%%%%
\bigskip
\section{Towards a Stability Spectrum Theorem} \label{s:stabspec}

In this section we begin work towards a stability spectrum theorem for
tame classes.  Other partial results in this direction appear in
\cite{BaKuVa}.

\begin{theorem}[Extension property for non-$\mu$-splitting types]
Let $\K$ be an AEC.  Suppose $N\prec_{\K}M'\prec_{\K}M$ are such that
$M'$ is universal over $N$ and all three models have cardinality $\mu$.
Then for every $p'\in\gaS(M')$ which does not $\mu$-split over $N$, there
exists $p\in\gaS(M)$ extending $p'$ such that $p$ does not $\mu$-split
over $N$.

\end{theorem}

\begin{proof}
Let $N,M',M$ be as in the statement of the theorem.
Suppose $p=\tp(a'/M')$  does not $\mu$-split over $N$.
Since $M'$ is universal over $N$, there exists $f:M'\rightarrow M$ with
$f\restriction N=\id_{N}$.  Notice that by monotonicity of
non-$\mu$-splitting, we have that $\tp(a'/f(M))$ does not $\mu$-split
over $N$.  Then, by invariance we have that also
$$(*)\quad\tp(f(a')/M)\text{ does not }
\mu\text{-split over }N.$$  By one more appeal to non-splitting, we can
conclude that $\tp(f(a')/M')=\tp(a'/M')$.  Otherwise these two types would
contradict $(*)$.
\end{proof}

\begin{theorem}[Uniqueness  of non-splitting
extensions]\label{unique non-split}
Let $\K$ be $\chi$-tame abstract elementary class and $\mu$ a cardinal
with
$\mu\geq\chi$.
 Let $M\in\K_{>\mu}$ and $M',N\in\K_\mu$ with
$N\prec_{\K}M'\prec_{\K}M$.
If $M'$ is universal over $N$,
then for every $p'\in \gaS(M')$ which does not $\mu$-split over $N$,
if there  are
$q,p\in
\gaS(M)$ such that
$q$ and $p$ both extend $p'$ and  do not $\mu$-split over $N$, then $p=q$.
\end{theorem}
\begin{proof}
 Suppose for the sake of contradiction that there
exists
$q\neq p\in\gaS(M)$ extending $p'\in\gaS(M')$ such that both $p$ and $q$
do not
$\mu$-split over
$N$.
Let $a, b$ be such that $p=\tp(a/M)$ and $q=\tp(b/M)$.
By tameness there exists $M^*\in\K_\mu$ such that
$M'\prec_{\K}M^*\prec_{\K}M$ and $p\restriction M^*\neq q\restriction
M^*$.

Since $M'$ is universal over $N$, there exists a $\prec_{\K}$-mapping
$f:M^*\rightarrow M'$ with $f\restriction N=\id_N$.
Since $p$ and $q$ do not $\mu$-split over $N$ we have
$$(*)_a\quad\tp(a/f(M^*))=\tp(f(a)/f(M^*))\text{ and}$$
$$(*)_b\quad\tp(b/f(M^*))=\tp(f(b)/f(M^*)).$$
On the other hand, since $p\restriction M^*\neq q\restriction M^*$, we
have that their images are also not equal
$$(*)\quad\tp(f(a)/f(M^*))\neq\tp(f(b)/f(M^*)).$$
Combining $(*)_a$, $(*)_b$ and $(*)$, we get
$$\tp(a/f(M^*))\neq\tp(b/f(M^*)).$$
Since $f(M^*)\prec_{\K}M'$, this inequality witnesses that
$$\tp(a/M')\neq\tp(b/M'),$$
contradicting our choice of $p$ and $q$ both extending $p'$.

\end{proof}

\begin{corollary}\label{few non splitting types} Let $\K$ be $\chi$-tame
abstract elementary  class and let $\mu$ be a cardinal with $\mu>\chi$.
%in Hypothesis \ref{tamehyp}
 Let $M\in\K_{>\mu}$, $N\prec_{\K}M$ with $M$ universal over $N$.
If $\K$ is Galois-stable in $\mu$ then
\[
| \{p\in \gaS(M)\;:\;  p \mbox{ does not }\mu\mbox{-split
over }N\}|\leq \mu.
\]

\end{corollary}

\begin{proof}
Since $M$ is universal over $N$, there is a $M'\in\K_\mu$
with
$N\prec_{\K}M'\prec_{\K}M$ and $M'$ universal over $N$. By $\mu$-stability
there are only
$\mu$-many types over
$M'$ and by Theorem \ref{unique non-split} each type over $N$ has at most
one non-$\mu$-splitting extension to $M$.
\end{proof}

\begin{corollary}[Partial stability spectrum]
Suppose that $\K$ is a tame abstract elementary class satisfying the
amalgamation property.  If $\K$ is Galois-stable in some $\mu>\Hanf(\K)$,
then
$\K$ is stable in every $\kappa$ with $\kappa^\mu=\kappa$.
\end{corollary}

\begin{proof}
Let $\kappa=\kappa^\mu$.
Consider $M\in\K_{\kappa}$ and let $M'$ be a $\mu$-saturated extension
of
$M$. Such an extension exists by stability in $\mu$ and our cardinal
arithmetic assumption.  Suppose that $\tp(a/M,\C)\neq\tp(b/M,\C)$, then
surely $\tp(a/M',\C)\neq\tp(b/M',\C)$. So the number of Galois-types over
$M$ is $\leq$ the number of Galois-types over $M'$.

Thus we may assume that $M$ is $\mu$-saturated.
Again by our assumption that $\kappa^\mu=\kappa$, we can enumerate
$\langle N_i\mid i<\kappa\rangle$ all the $\prec_{\K}$-submodels of $M$ of
cardinality $\mu$.  Since $M$ is $\mu$-saturated, for each one of the
$N_i$ there is $M'_i$ of cardinality $\mu$ such that
$N_i\prec_{\K}M'_i\prec_{\K}M$ with $M'_i$ universal over $N_i$.  By Fact
\ref{bounded kapp} for every
$p\in\gaS(M)$, there exists
$N_i\prec_{\K}M$ of cardinality $\mu$ such that $p$ does not
$\mu$-split over $N_i$.
Thus $$|\gaS(M)|=|\Union_{i<\kappa}\{p\in\gaS(M)\mid p\text{ does not
}\mu\text{-split over }N_i\}|.$$

By  monotonicity of non-splitting we have that for each $p\in\gaS(M)$
that does not $\mu$-split over $N_i$,  also $p\restriction M'_i$ does not
$\mu$-split over $N_i$.
By
Corollary
\ref{few non splitting types}, for each $i<\kappa$,
$$|\{p\in\gaS(M)\mid p \text{ does not }\mu\text{-split over }N_i\}|$$
$$\leq
|\{p\in\gaS(M'_i)\mid p \text{ does not }\mu\text{-split over }N_i\}|.
$$
Now by stability in $\mu$ and our assumption that $\kappa^\mu=\kappa$, we
can conclude that $|\gaS(M)|\leq\kappa$.

\end{proof}

%%%%%%%%%%%%%%%%%%%%%%%%%%%%%%%%%%%%%%%%%%%%%%%%%%%%%%%%%%%%%%%%%%
\bigskip
\section{Exercise on Dividing} \label{s:dividing}

With the existence of Morley sequences a natural extension is to study the
following dependence relation to determine whether or not it satisfies
properties such as transitivity, symmetry or extension.  Here we derive
the existence property for types over saturated models.

\begin{definition}\label{divides}
Let $p\in \gaS(M)$ and $N\prec_{\K}M$.  We say that \emph{$p$ divides over
$N$} \index{divides}\index{Galois-type!divides} iff there are $a\in
|M|\backslash|N|$ and a Morley sequence,
$\{a_n\mid n<\omega\}$ for the $\tp(a/N,M)$ such that for
some collection\\
$\{f_n\in Aut_M\C\mid n<\omega\}$ with $f_n(a)=a_n$ we have
$$\{f_n(p)\mid n<\omega\}\text{ is inconsistent}.$$
I.e. there is no common realization to the above mentioned $\omega$-many
types.
\end{definition}

\begin{theorem}[Existence]\index{divides!existence}\label{divides!existence}
Suppose that $\K$ is Galois-stable in $\mu$ and $\chi$-tame for some
$\chi<\mu$.  For every $\mu$-model homogeneous model
$M\in\K_{\geq\mu}$ and ever $p\in
\gaS(M)$, there  exists $N\prec_{\K}M$ of cardinality $\mu$ such
that $p$ does not divide over $N$.
\end{theorem}

\begin{proof}
Suppose that $p$ and $M$ form a counter-example.  We will find $N_i,
N_i^1, N_i^2$ and $ h_i$ for
$i<\mu$ satisfying:
\begin{enumerate}
\item $\langle N_i\in\K_\mu\mid i\leq\mu\rangle$ is a
$\prec_{\K}$-increasing and continuous sequence of models;
\item $N_i\prec_{\K}N_i^l\prec_{\K}N_{i+1}$ for $i<\mu$ and $l=1,2$;
\item for $i<\mu$, $h_i:N_i^1\cong N_i^2$ and $h_i\restriction
N_i=\id_{N_i}$ and
\item $p\restriction N_i^2\neq h_i(p\restriction N_i^1)$.
\end{enumerate}

Suppose that $N_i$ has been defined.  Since $p$ divides over every
substructure of cardinality $\mu$, we may find $\bar a$,
$\{\bar a_n\mid n<\omega\}$ and $\{f_n\mid n<\omega\}$ witnessing that
$p$ divides over $N_i$.  Namely, we have that $\{f_n(p)\mid n<\omega\}$
is inconsistent.  Let $n<\omega$ be such that $f_0(p)\neq f_n(p)$.
Then $p\neq f_0^{-1}\circ f_n(p)$.  By $\chi$-tameness, we can find
$N^*\prec_{\K}M$ of cardinality $\mu$ containing $N$ such that
$p\restriction N^*\neq (f_0^{-1}\circ f_n(p))\restriction N^*$.  WLOG
$f_0^{-1}\circ f_n\in Aut_NN^*$.

Let $h_i:=f_0^{-1}\circ f_n$, $N_i^1:=N^*$ and $N_i^2:=N^*$.  Choose
$N_{i+1}\prec_{\K}M$ to be an extension of $N^*$ of cardinality $\mu$.

Now we can use   $\l  N_i, N_i^1, N_i^2, h_i\mid i<\mu\r$ to contradict
Galois-stability in $\mu$ by repeating the argument from the proof of
Fact
\ref{bounded kapp}, (starting at the second paragraph) we construct an
increasing chain of models and a tree of $\K$-embedding which by Claim
\ref{many types claim} gives many Galois-types contradicting stability
in $\mu$.
\end{proof}

\begin{remark}
In a previous version of this paper, we stated for the existence property
for arbitrary models $M\in\K_{\geq\mu}$.  We do not know how to remove
the assumption of $\mu$-model homogeneity and we suspect that this is not
possible with only the assumption of $\mu$-stability and tameness.
\end{remark}

%\begin{remark}
% Buechler and Lessmann in  \cite{BuLe} introduced the notion of
%simplicity for the context of a
%given sequentially homogenous model.  The definition is that dividing (of
%first-order formulas)  satisfies two requirements:
%\begin{enumerate}
%\item
%Existence and
%
%\item
%Extension property.
%
%\end{enumerate}
%
%Missing from the preliminary draft was a proof  that stability
%implies simplicity for homogenous models.  After they circulated a
%preprint, Shelah discovered an example of a homogeneous models which is
%stable but not simple!  The example had the existence property but failed
%to have the extension property.  A biproduct of Theorem
%\ref{divides!existence} is that
% existence always
%follows from stability (even for a much more general contexts than
%homogeneous models).  Notice that in homogeneous model theory we
%essentially require that every set is an amalgamation base (not only
%models like in this paper).
%Shelah's example show that Galois-stability is insufficient to imply the
%extension property for Galois-types.
%
%
%\end{remark}

%%%%%%%%%%%%%%%%%%%%%%%%%%%%%%%%%%%%%%%%%%%%%%%%%%%%%%%%%%%%%%%%%%
%%%%%%%%%%%%%%%%%%%%%%


\begin{thebibliography}{10}
\bibitem[AlGr]{AlGr} M. Albert and R. Grossberg,
\newblock Rich models,
\newblock
Journal of Symbolic Logic, 55, (1990) 1292--1298.

\bibitem[Ba]{Ba}
John Baldwin.
\newblock The Hrushovski construction and quasiminimal excellent classes.
\newblock In preparation.

\bibitem[BaKuVa]{BaKuVa} J. Baldwin, D. Kueker and M. VanDieren,
\newblock Upward Stability Transfer Theorem for Tame Abstract Elementary
Classes, (8 pages), Preprint available
\texttt{http://www.math.lsa.umich.edu/\~{}mvd/home.html}.

\bibitem[BaSh]{BaSh}
John Baldwin  and S. Shelah.
\newblock
Examples of non-locality. in preparation.

%\bibitem[Be]{Be}
%A. Berenstein,  Some generalizations of first order
% tools to homogeneous models,  To appear in \emph{J.S.L.},  Preprint
%available at
%\texttt{http://www.math.uiuc.edu/\~{}aberenst/research.html}


%\bibitem[BeBu]{BeBu}
%A. Berenstein and S. Buechler,  A study of independence in
%strongly homogeneous expansions of Hilbert spaces,  Preprint available at
%\texttt{http://www.math.uiuc.edu/\~{}aberenst/research.html}


%\bibitem[BuLe]{BuLe}
%S. Buechler and O. Lessmann,  Simple homogeneous models,
%\emph{J. Amer. Math. Soc.}, {\bf 16} (2003), no. 1, 91--121.

%\bibitem[BY]{BY}
%I. Ben-Yaacov,  Positive model theory and compact
%abstract theories, \emph{J. Math. Log.}, {\bf 3} (2003), no. 1, 85--118.

%\bibitem[CK]{CK}
%C. ~C. Chang and H.~Jerome  Keisler,
%\newblock \textbf{ Model Theory},
%\newblock North-Holland Publishing Co., Amsterdam, 1990.


%\bibitem[Fr]{Fr}  R. Fra\"{i}ss\'{e},
%Sur quelques classifications des syst\'{e}mes de relations,
% {\em Publ. Sci. Univ. Alger. S\'{e}r A} {\bf 1} (1954), 35--182.



\bibitem[Gr1]{Gr1}
R. Grossberg,
Classification theory for non-elementary classes,
{\bf Logic and Algebra}, ed. Yi Zhang,
{\em Contemporary Mathematics}, {\bf 302}, (2002) AMS,  pp. 165--204.



\bibitem[Gr2]{Gr2}
R. Grossberg,
\newblock {\bf A Course in Model Theory,}
\newblock Book in Preparation, Available at
\texttt{http://www.math.cmu.edu/\~rami/home.html}.



\bibitem[GrLe1]{GrLe1}
R. Grossberg and O. Lessmann,
The local order property in non
elementary      classes,
{\em Arch Math Logic} {\bf 39} (2000) 6, 439--457.



\bibitem[GrLe2]{GrLe2}
R. Grossberg and O. Lessmann,
Shelah's stability spectrum and
homogeneity     spectrum in finite diagrams,
 {\em  Archive for
Mathematical Logic}, {\bf 41}, (2002) 1, 1--31.

\bibitem[GrKo]{GrKo}
R. Grossberg and A. Kolesnikov,
Excellent  abstract
elementary classes are tame, (12 pages), preprint.


%\bibitem[GrSh]{GrSh}
%R. Grossberg and S. Shelah,
%\newblock On universal locally finite groups
%\newblock \emph{ Israel J. of Math. }, {\textbf 44}, 1983, 289--302.


\bibitem[GrVa1]{GrVa1} R. Grossberg and M. VanDieren,
Categoricity in Tame Abstract Elementary Classes, (23 pages), Preprint
available at
\texttt{http://www.math.lsa.umich.edu/\~{}mvd/home.html}.

\bibitem[GrVa2]{GrVa2} R. Grossberg and M. VanDieren,
Categoricity from one successor cardinal               in Tame Abstract
Elementary Classes,  (17 pages),  Preprint available at
\texttt{http://www.math.lsa.umich.edu/\~{}mvd/home.html}.

\bibitem[GVV]{GVV} Rami Grossberg, Monica VanDieren and Andr\'es
Villaveces,  Limit Models in Classes with Amalgamation, (18 pages),
Preprint available at
\texttt{http://www.math.lsa.umich.edu/\~{}mvd/home.html}.

\bibitem[Ha]{Ha}
V. Harnik,
\newblock On existence of saturated models
of stable theories,
\newblock {\em Proc. of AMS}, {\bf 52}, 361--367, 1975.

%\bibitem[Ho]{Ho} Wilfrid Hodges.
%{\bf Model Theory},
%Cambridge University Press, 772 pages,  1993.

%\bibitem[Hy]{Hy}
%T. Hyttinen,
%\newblock On nonstructure of elementary submodels of a stable
%homogeneous
%structure,
%\newblock {\em Fundamenta Mathematicae}, {\bf 156}, 1998, 167-182.

%\bibitem[HySh]{HySh}
%T. Hyttinen and S. Shelah,
%Strong splitting in stable homogeneous
%models,
%\emph{APAL}, {\bf 103}, 2000, 201--228.


%\bibitem[Jo1]{Jo1} B. J\'{o}nsson,
%\newblock Universal relational systems,
%\newblock {\em Math. Scand.}   {\bf 4} (1956) 193--208.


%\bibitem[Jo2]{Jo2}
%B. J\'{o}nsson,
%\newblock Homogeneous universal systems,
%\newblock {\em Math. Scand.}, {\bf 8}, 1960,  137--142.

%\bibitem[Ke1]{Ke1} H. J. Keisler,
%\newblock $L_{\omega_1,\omega}(\mathbf{Q})$
%\newblock \emph{ Ann of Math Logic}, {\textbf 1}, 1969.

%\bibitem[Ke2]{Ke2} H. J. Keisler,
%\newblock {\bf Model Theory for Infinitary Logic},
%\newblock North-Holland 1971.


%\bibitem[LaPo]{LaPo}
%Daniel Lascar and Bruno Poizat.
%\newblock An introduction to forking,
%{\em Journal of Symbolic Logic}, {\bf 44}, 1979, 330--350.

%\bibitem[Le]{Le}
%O. Lessmann,
%Pregeometries in finite diagrams.
%\relax {\em Ann. Pure Appl. Logic}, {\bf 106},  no. 1-3, 49--83, 2000.

%\bibitem[MaSh]{MaSh} A. Macintyre and S. Shelah, Uncountable
%universal locally finite groups.  \emph{J. of  Algebra} {\bf  43} (1976)
%168--175.

\bibitem[Ma]{Ma} L. Marcus,
A prime minimal model with an infinite set of indiscernibles,
{\em Israel Journal of Mathematics},
{\bf 11}, (1972), 180-183.

\bibitem[MaSh]{Sh285}
Michael Makkai and Saharon Shelah.
\relax  {Categoricity of theories in $L_ {\kappa\omega},$ with $\kappa$ a
     compact cardinal}.
\relax {\em {Annals of Pure and Applied Logic}}, {\bf 47}:41--97, 1990.


%\bibitem[Rob56]{Rob56} Abraham Robinson.
%    A result on consistency and its application to the theory of
%              definition,
%  \emph{Nederl. Akad. Wetensch. Proc. Ser. A. {\bf 59}},  = Indag. Math.,
%  {\bf 18}}, {47--58},
%    {1956}.

\bibitem[Shc]{Shc}
S. Shelah,
\newblock {\bf Classification Theory and the
Number of Non-isomorphic Models} $2^{nd}$ edition,
\newblock {North Holland Amsterdam}, 1990.


\bibitem[Sh3]{Sh3}
S. Shelah,
\newblock Finite diagrams stable in power,
\newblock {\em Ann. Math. Logic}, {\bf 2}, 69--118, 1970/1971.


%\bibitem[Sh 54]{Sh54}
%S. Shelah,
%\newblock The lazy model-theoretician's guide to stability,
%\newblock {\em Logique et
%Analyse} {\bf 18} (1975) 241-308.


\bibitem[Sh 87a]{Sh 87a}
S. Shelah,
\newblock Classification theory for nonelementary classes, I,
The number of uncountable models of $\psi \in L_{\omega _{1},\omega }$,
Part A,
\newblock {\em Israel J. Math.}, {\bf 46}:212--240, 1983.

\bibitem[Sh 87b]{Sh 87b}
S. Shelah,
\newblock Classification theory for nonelementary classes, I,
The number of uncountable models of $\psi \in L_{\omega _{1},\omega }$,
Part B,
\newblock {\em Israel J. Math.}, {\bf 46}:241--273, 1983.



\bibitem[Sh 88]{Sh 88}
S. Shelah,
\newblock Classification of nonelementary classes, II,
Abstract Elementary Classes.  In
{\bf Classification Theory (Chicago IL 1985)}, volume 1292 of
\emph{Lecture Notes in Mathematics}, pages 419--497.  Springer, BSerlin, 
1987.
Proceedings of the USA-Israel Conference on Classification Theory, Chicago,
December 1985; ed. Baldwin, J.T.


\bibitem[Sh 300]{Sh 300}
S. Shelah,
\newblock Universal Classes.
\newblock In \emph{Classification Theory} of
\emph{Lecture Notes in Mathematics},
{\bf 1292}, 264-418.
Springer-Berlin, 1987.



\bibitem[Sh 394]{Sh 394}
\newblock S. Shelah,
\newblock Categoricity of abstract classes with amalgamation,
\newblock \emph{Annals of Pure and Applied Logic}, {\bf 98}(1-3),
pages 141--187, 1999.

\bibitem[Sh 576]{Sh 576}
S. Shelah,
\newblock Categoricity of an abstract elementary class in two
successive cardinals,
\newblock \emph{Israel J. of Math},  {\bf 126}, (2001), 29--128.

\bibitem[Sh 600]{Sh 600}
S. Shelah,
\newblock Categoricity in abstract elementary classes: going up inductive
step,
\newblock Preprint,  (100 pages).

\bibitem[Sh 702]{Sh702}
S. Shelah,
\newblock
On what I do not understand (and have something to say), model theory,
\emph{Math. Japon.}, {\bf 51}, (2000), no. 2, 329--377.

\bibitem[Sh 705]{Sh705}
S. Shelah,
\newblock
Toward classification theory of
good $\lambda$ frames and abstract elementary classes,  In preparation.


\bibitem[ShVi]{ShVi 635}
S. Shelah and A. Villaveces,
\newblock Categoricity in abstract elementary classes with no maximal
models,
\newblock \emph{Annals of Pure and Applied Logic}, 1999.

\bibitem[Va]{Va}
M. VanDieren,
\newblock  Categoricity in abstract elementary classes with no maximal
models,
\newblock (57 pages),  accepted, subject to revisions by APAL.
\newblock Preprint available at
\texttt{www.math.lsa.umich.edu/\~{}mvd/home.html}


\bibitem[ViZa]{ViZa}
Andr\'es Villaveces and Pedro Zambrano.
\newblock Hrushovski constructions and tame abstract elementary classes.
\newblock in preparation.


\bibitem[Zi]{Zi}
B. Zilber,
\relax  {Analytic and pseudo-analytic structures,}
Preprint available at
\texttt{www.maths.ox.ac.uk/\~{}zilber}




\end{thebibliography}
\end{document}